\documentclass[3pt]{elsarticle}
\pdfoutput=1
\usepackage{amsthm}
\usepackage[intlimits]{amsmath}
\usepackage{amssymb}
\usepackage{amsthm}
\usepackage{amsfonts}
\usepackage{bm}
\usepackage{mathdots}

\usepackage[utf8]{inputenc} % for UTF-8 encoding
\usepackage[T1]{fontenc}
\usepackage{times}
\usepackage{mathptmx}  % times roman, including math
\DeclareSymbolFont{cmletters}{OML}{cmr}{m}{n}
\DeclareMathAlphabet{\mathcal}{OMS}{cmsy}{m}{n} % recover 'mathcal'

\journal{Computer Methods in Applied Mechanics and Engineering}
\ifdefined\copyReadyOn
\usepackage{ecrc}
\volume{00}
\firstpage{1}
\runauth{G. Beer, Ch. Duenser, E. Ruocco, V. Mallardo}
\journalname{\journal}
\jid{CMAME}
\jnltitlelogo{CMAME}
\CopyrightLine{2015}{Published by Elsevier Ltd.}
\fi

\usepackage{tikz}
\usepackage{tikz-3dplot}
\usepackage{pgfplots}
\usepackage{pgfplotstable}
\usepackage{tikzscale}
\usepackage{standalone} % for standalone compilation of tikz graphics
\usepackage{hf-tikz}

\usepackage{forloop}
\usepackage{arrayjobx}

\usepackage{algorithmic}
\usepackage{algorithm}

\usepackage{booktabs} % for tables: toprule, midrule, bottomrule
\usepackage{calc}
\usepackage{xcolor,pdfcomment,soul} % add pdf comments
\usepackage{xspace} % horizontal space after makro words
\usepackage{multirow}

\usepackage{scrpage2}
\usepackage{setspace}

\usepackage[labelformat=simple]{subcaption}   	%% 1(a)
\usepackage{pdfpages} % to include external pdf (i.e. eidesstattliche erlärung)

\usepackage{tabto} % for tabbing within itemize

\usepackage{makeidx}
\makeindex

\usepackage{colortbl}

\usepackage{enumerate}

\usepackage{overpic}

\usepackage{listings} %% enviroment for code formating 

\usepackage[intlimits]{amsmath}
\usepackage{amssymb} 
\usepackage{amsthm}
\usepackage{amsfonts,mathrsfs}
\usepackage{bm}
\usepackage{mathtools} % for \coloneqq and \shortintertext

\usepackage{siunitx} %% units

\usepackage{xfrac}
\usepackage{listings}
\usepackage{color}
\definecolor{dkgreen}{rgb}{0,0.6,0}
\definecolor{gray}{rgb}{0.5,0.5,0.5}
\definecolor{mauve}{rgb}{0.58,0,0.82}

\usepackage[utf8]{inputenc} % for UTF-8 encoding
\usepackage[german,english]{babel}  
\usepackage[T1]{fontenc} % größerer Zeichenvorat z.B. für Umlaute
\usepackage{lmodern}

\usepackage{geometry}
\usepackage{overpic}
\graphicspath{{pics/}{tmp/}}

\newboolean{inGray}
\setboolean{inGray}{false}
\newcommand{\addtoindex}[2][]{% key=#1,text=#2
    \ifthenelse { \equal{#1}{} }
    {#2\index{#2}\xspace}%
    {#2\index{#1}\xspace}%  
}

\newcommand{\myVec}[1]{\mathbf{#1}}
\newcommand{\myVecGreek}[1]{\boldsymbol{#1}}

\newcommand\tens[2]{\mathsf{#1}} % 
\newcommand{\trans}{\text{T}} % transposed matrix/vector
\newcommand{\NURBS}{R} 			% symbol for NURBS
\newcommand\uu{\xi} 			% Local NURBS coordinates
\newcommand\vv{\eta} 			% Local NURBS coordinates
\newcommand{\domain}{\Omega}

\newcommand{\boundary}{\Gamma}

\newcommand\primary{u} 			% primary variable - displacements in elasticity
\newcommand\dual{t} 			% dual variable - tractions in elasticity
\providecommand\url[1]{\emph{#1}}

\newcommand\fund[1]{\tens{#1}{2}}

\newcommand\pt[1]{\boldsymbol{#1}}
\newcommand\sourcept{\tilde{\pt{x}}}
\newcommand\fieldpt{\hat{\pt{x}}}

{
  \begin{table}[#1]%[htb!]
    \def\mycap{#2}
    \def\mylabel{#3}
    \centering
    \begin{tabular}{#4}
      \toprule
}
{
  \bottomrule
  \end{tabular}
  \caption{\mycap}
  \label{\mylabel}
  \end{table}
}

{
  \begin{algorithm}
    \caption{#1}
    \label{#2}
    \begin{algorithmic}[1]
}   
{
\end{algorithmic}
\end{algorithm}
}

\newtheoremstyle{myremark}% ⟨name⟩ 
{3pt}% Space above
{3pt}% Space below
{}% Body font
{}% Indent amount
{\itshape}% Theorem head font
{:}% Punctuation after theorem head
{.5em}% Space after theorem head
{}% Theorem head spec (can be left empty, meaning ‘normal’)
\theoremstyle{myremark}
\newtheorem*{remark}{Remark}

\newcommand{\myfigref}[1]{Figure~\ref{#1}}

\newcounter{footnoteNumber} % initialization of counter with 0

\DeclareCaptionFormat{listing}{#1#2#3}
\newcommand{\myalignatsinglelabel}[3]%
{%
    \begin{equation}
        \label{#2}
        \begin{alignedat}{#1} 
            #3
        \end{alignedat}
    \end{equation}
}

\newcommand{\myalignat}[2]%
{%
  \begin{alignat}{#1} 
    #2
   \end{alignat}
}

\usepackage{pgf,pgffor,pgfmath}
\pgfplotsset{compat=1.10}

\usetikzlibrary{external}
\usetikzlibrary{positioning,calc,intersections,arrows,3d}
\usetikzlibrary{shapes,snakes,fit,spy}
\usepgfplotslibrary{patchplots}
\usetikzlibrary{matrix}

\usepgfplotslibrary{groupplots}
\usepgfplotslibrary{fillbetween}

\def\pgfplotfontsizetitle{\small}
\def\pgfplotfontsizelegend{\small}
\def\pgfplotfontsize{\small}
\def\pgfplotfontsizetiny{\scriptsize}

\def\tikzfontsizetiny{\scriptsize}

\pgfplotsset{
  mystyle/.style ={%
    grid = major,
    every tick label/.append style={font=\pgfplotfontsizetiny},
    every axis label/.append style={font=\pgfplotfontsize},
    legend style={font=\pgfplotfontsizelegend},
    label style={font=\pgfplotfontsize},
    title style={font=\pgfplotfontsizetitle},
    /pgf/number format/set thousands separator = {}, % 1,000 -> 1000
  }
}% 

\makeatletter
\pgfplotsset{
    myIgnoreRowModulo2/.style args={#1}{%
        /pgfplots/x filter/.code={%
        \let\xValue\pgfmathresult % store x value of current table index
        \pgfmathparse{int(mod(int(\coordindex),int(2))} \pgfmathresult % compute modulo
        \ifnum#1=\pgfmathresult
            \def\pgfmathresult{} % ignore table entry
        \else
            \pgfmathparse{\xValue} \pgfmathresult % restore x value to \pgfmathresult macro
        \fi
        }
    } 
}
\makeatother

\colorlet{drawblue}      {blue!80!white}
\colorlet{drawred}       {red!80!white}
\colorlet{drawgray}      {gray}
\definecolor{drawgreen}  {RGB}{44,162,95}
\colorlet{drawpurple}    {purple}
\colorlet{draworange}    {orange}
\colorlet{drawlime}      {lime!80!black}
\colorlet{drawartichoke} {yellow!60!black}
\colorlet{TUGgray}{black!15}
\definecolor{TUGred}{RGB}{247,1,70}
\definecolor{IFBblue}{RGB}{51,112,169}

\definecolor{basisColor1}{RGB}{59,76,192}
\definecolor{basisColor2}{RGB}{87,117,225}
\definecolor{basisColor3}{RGB}{119,154,247}
\definecolor{basisColor4}{RGB}{152,185,255}
\definecolor{basisColor5}{RGB}{184,208,249}
\definecolor{basisColor6}{RGB}{195,209,230}	% changed for printing
\definecolor{basisColor7}{RGB}{247,200,190}	% changed for printing
\definecolor{basisColor8}{RGB}{247,187,160}
\definecolor{basisColor9}{RGB}{244,154,123}
\definecolor{basisColor10}{RGB}{229,112,88}
\definecolor{basisColor11}{RGB}{203,62,56}
\definecolor{basisColor12}{RGB}{180,4,38}

\definecolor{basisColor8sw}{RGB} {189,189,189}
\definecolor{basisColor9sw}{RGB} {150,150,150}
\definecolor{basisColor10sw}{RGB}{115,115,115}
\definecolor{basisColor11sw}{RGB}{91,91,91}
\definecolor{basisColor12sw}{RGB}{37,37,37}

\colorlet{myblue}    {blue}
\colorlet{myred}     {red}
\colorlet{mygreen}   {drawgreen}
\colorlet{mypurple}  {purple}
\colorlet{myorange}  {orange}

\ifthenelse{\boolean{inGray}}{ 
\pgfplotscreateplotcyclelist{mycyclelist}{
  basisColor12sw,semithick, every mark/.append style={solid, fill=.},mark=*           \\%
  basisColor11sw,semithick, every mark/.append style={solid, fill=.},mark=square*     \\%
  basisColor10sw,semithick, every mark/.append style={solid, fill=.},mark=triangle*   \\%
  basisColor9sw , semithick,every mark/.append style={solid, fill=.},mark=diamond*    \\%
  basisColor8sw , semithick,every mark/.append style={solid, fill=.},mark=otimes    \\%
}
}{

\pgfplotscreateplotcyclelist{mycyclelist}{
  basisColor12,semithick,   every mark/.append style={solid, fill=.},mark=*           \\%
  basisColor11,semithick,   every mark/.append style={solid, fill=.},mark=square*     \\%
  basisColor10,semithick,   every mark/.append style={solid, fill=.},mark=triangle*   \\%
  basisColor9, semithick,   every mark/.append style={solid, fill=.},mark=diamond*    \\%
  basisColor8, semithick,   every mark/.append style={solid, fill=.},mark=otimes    \\%
}
}

\ifthenelse{\boolean{inGray}}{ 
\pgfplotscreateplotcyclelist{mycyclelistcompare}{
  basisColor12sw, semithick, loosely dashdotdotted     \\%
  basisColor11sw, semithick, loosely dashed      		\\%
  basisColor10sw, semithick, dashed \\%loosely dashdotted        \\%
  basisColor9sw,  semithick, solid               	\\%
  basisColor12sw, only marks,   every mark/.append style={solid, fill=.},mark=otimes    \\%
  basisColor11sw, only marks,   every mark/.append style={solid, fill=.},mark=triangle*   \\%
  basisColor10sw, only marks,   every mark/.append style={solid, fill=.},mark=square*     \\%
  basisColor9sw,  only marks,   every mark/.append style={solid, fill=.},mark=*          \\%
}

}{
\pgfplotscreateplotcyclelist{mycyclelistcompare}{
  basisColor12, semithick, loosely dashdotdotted     \\%
  basisColor11, semithick, loosely dashed      		\\%
  basisColor10, semithick, dashed \\%loosely dashdotted        \\%
  basisColor9,  semithick, solid               	\\%
  basisColor12, only marks,   every mark/.append style={solid, fill=.},mark=otimes    \\%
  basisColor11, only marks,   every mark/.append style={solid, fill=.},mark=triangle*   \\%
  basisColor10, only marks,   every mark/.append style={solid, fill=.},mark=square*     \\%
  basisColor9,  only marks,   every mark/.append style={solid, fill=.},mark=*          \\%
}
}

\ifthenelse{\boolean{inGray}}{ 
\tikzset{mycyclelistcompareReferenceA/.style={basisColor12sw,solid}}
\tikzset{mycyclelistcompareTestA/.style={basisColor12sw,only marks,mark=otimes}}
}{
\tikzset{mycyclelistcompareReferenceA/.style={basisColor12,solid}}
\tikzset{mycyclelistcompareTestA/.style={basisColor12,only marks,mark=otimes}}
}
\tikzset{helpline/.style={thin,dashed}}
\tikzset{labelline/.style={thin}}
\tikzset{referencePath/.style={dotted,very thick}}
\tikzset{helparrow/.style={thin,arrows={-latex}}}
\tikzset{axis/.style={thin,arrows={->}}}
\tikzset{force/.style={thick,arrows={->}}}
\tikzset{forceInverse/.style={thick,arrows={<-}}}
\tikzset{Gamma/.style={ultra thick}}
\tikzset{controlPoly/.style={draw=black}}
\tikzset{GammaFill/.style={fill=lightgray,fill opacity=0.5}}
\tikzset{colorDiri/.style={drawgreen}}
\tikzset{GammaFillDiri/.style={fill=drawgreen,fill opacity=0.5}}
\tikzset{initialgrid/.style={thin,gray}}
\tikzset{addgridline/.style={dashed,gray}}
\tikzset{trimmingcurve/.style={thick}}
\tikzset{boundingbox/.style={thick, dotted}}
\tikzset{parameterSpace/.style={ }}
\tikzset{basisfunction/.style={very thick,smooth}}
\tikzset{bspline/.style={very thick,smooth,red}}
\tikzset{intersectioncurve/.style={dashed,thick}}
\tikzset{integrationRegionEdge/.style={dashed}}
\tikzset{pointer/.style={arrows={-latex}}}

\tikzstyle{anode}= [circle, inner sep=1.3pt, draw, fill=black]
\tikzstyle{gausspoint}=[shape=circle,draw=black,fill=black,inner sep=1.1pt]
\tikzstyle{controlPoint}=[shape=circle,draw=black,fill=black,thin,inner sep=0pt,minimum size=1.5mm]
\tikzstyle{abscissaPoint}=[shape=circle,draw=black,fill=white,thin,inner sep=0pt,minimum size=1.5mm]
\tikzstyle{anchorPoint}=[shape=circle,draw=black,fill=black,thin,inner sep=0pt,minimum size=1.5mm]
\tikzstyle{anchorPointDeg}=[shape=circle,draw=black,fill=TUGred,thin,inner sep=0pt,minimum size=1.5mm]
\tikzstyle{anchorPointDegD}=[shape=cross out,thick,draw=black,inner sep=0pt,minimum size=1.5mm]
\tikzstyle{trimmingIntersectionPoint}=[shape=cross out,thick,draw=black,inner sep=0pt,minimum size=1.5mm]

\tikzset{%
  highlight/.style={rectangle,rounded corners,fill=red!60,draw,fill opacity=0.125,thick,inner sep=0pt}
}

\usepackage{colortbl}

\def\trianglecolor{black}
\newcommand{\upperSlopeTriangle}[4] 	% input: #1 slope #2 vertical shift  #3 min x-value #4 max x-value
{
	\def\trianglecolor{black}
	\addplot[forget plot, domain=#3:#4,color=\trianglecolor,samples=2]{  #2 / (x^#1) } node (A1) [pos=1] {}; 
	\addplot[forget plot, domain=#3:#4,color=\trianglecolor,samples=2]{   #2  / (#3^#1)} node (A2) [pos=1] {} node [anchor=south,pos=0.5,black] {\tikzfontsizetiny $1$};
	\draw[color=\trianglecolor] (A1.center) -- (A2.center) node [anchor=west,pos=0.5,black] {\tikzfontsizetiny #1};
}

\newcommand{\lowerSlopeTriangle}[4] 	% input: #1 slope #2 vertical shift  #3 min x-value #4 max x-value
{
	\def\trianglecolor{black}
	\addplot[forget plot, domain=#3:#4,color=\trianglecolor,samples=2]{  #2 / (x^#1) } node (A1) [pos=0] {}; 
	\addplot[forget plot, domain=#3:#4,color=\trianglecolor,samples=2]{   #2  / (#4^#1)} node (A2) [pos=0] {} node [anchor=north,pos=0.5,black] {\tikzfontsizetiny $1$};
	\draw[color=\trianglecolor] (A1.center) -- (A2.center) node [anchor=east,pos=0.5,black] {\tikzfontsizetiny #1};
}

\newcommand{\myaddgraphic}[5]
{
 \node[anchor=south west,inner sep=0] (image) {\phantom{\includegraphics[#2]{#1}}};
  \begin{scope}[x={(image.south east)},y={(image.north west)}]
      
      \begin{scope}
          
          #5
          
          \node[anchor=south west,inner sep=0] {\includegraphics[#2]{#1}};
      \end{scope} 
      
      #4
      
      \pgfmathparse{int(#3)} \let\gridIndicator\pgfmathresult
      \ifthenelse{ \gridIndicator = 1 }
      {
          \draw[help lines,xstep=.1,ystep=.1] (0,0) grid (1.001,1.001);
          \foreach \x in {1,...,9} { \node [anchor=north] at (\x/10,0) {\x};}
          \foreach \y in {1,...,9} { \node [anchor=east] at (0,\y/10) {\y};}
      }{}
      
  \end{scope}    
}

\tikzstyle{reverseclip}=[insert path={(current page.north east) --
    (current page.south east) --
    (current page.south west) --
    (current page.north west) --
    (current page.north east)}
]

\newcounter{itR}
\newcommand{\bsplinevalue}[5] % input: #1 knot vector #2 order,  #3 intrinsic coordinate #4 span index #5 output name
{                    				
	\newarray\vKnots
	\newarray\vN
	\newarray\vNumeratorL
	\newarray\vNumeratorR
	\newarray\vSave

	\readarray{vKnots}{#1}
	\readarray{vN}{1}
	\readarray{vSave}{0}
	\readarray{vNumeratorL}{0}
	\readarray{vNumeratorR}{0}

	\foreach \j in {1,...,#2}
	{        
		\pgfmathparse{ int(#4+\j+1) } \checkvKnots(\pgfmathresult)
		\pgfmathsetmacro{\numR}{\cachedata-#3}
		
		\pgfmathparse{ int(#4-\j+2) } \checkvKnots(\pgfmathresult)
		\pgfmathsetmacro{\numL}{#3-\cachedata} 					

		\expandarrayelementtrue
		\pgfmathparse{ int(\j+1) }
		\vNumeratorL(\pgfmathresult)={\numL}
		\vNumeratorR(\pgfmathresult)={\numR}
		
		\forloop[1]{itR}{0}{\value{itR} < \j }
		{
			\pgfmathparse{ int(\theitR+1) } \let\tS\pgfmathresult  	
			\checkvSave(\tS)							
			\pgfmathsetmacro{\save}{\cachedata}  
		
			\pgfmathparse{int(\j-\theitR+1)} \checkvNumeratorL(\pgfmathresult)
			\pgfmathsetmacro{\tmpL}{\cachedata} 	
		    
			\pgfmathparse{ int(\theitR+2) } \checkvNumeratorR(\pgfmathresult)
			\pgfmathsetmacro{\tmpR}{\cachedata} 	       
	
			\pgfmathparse{ int(\theitR+1) } \checkvN(\pgfmathresult)
			\pgfmathparse{ \cachedata / (\tmpL+\tmpR) } \let\tmp\pgfmathresult
			
			\pgfmathparse{ \save + \tmpR * \tmp } 
			\vN(\tS)={\pgfmathresult}

			\pgfmathparse{ \tmpL * \tmp } \let\tmpsave\pgfmathresult
			\pgfmathparse{ int(\theitR+2) } \let\tS\pgfmathresult      
			\vSave(\tS)={\tmpsave}
		}
		
		\pgfmathparse{ int(\j+1) } \checkvSave(\pgfmathresult )
		\vN(\tS)={\cachedata}
	}

	\pgfmathparse{int( #2+1) } \let\lastIndex\pgfmathresult
	\checkvN(1) \pgfmathsetmacro{\first}{\cachedata}  
	\foreach \i [remember=\a as \lasta (initially \first)] in {2,...,\lastIndex}
	{
		\checkvN(\i) \def\a{\lasta,\cachedata} 
		\ifthenelse{\i=\lastIndex}{ \xdef#5{\a} }{}
	}
    
	\foreach \i in {1,...,\lastIndex}
	{
	    \clrarray{vNumeratorR}(\i)
	    \clrarray{vNumeratorL}(\i)
	    \clrarray{vN}(\i)
	    \clrarray{vSave}(\i)
	}
	\delarray\vN
	\delarray\vNumeratorL
	\delarray\vNumeratorR
	\delarray\vSave
}

\newcounter{countvalues}
\newcounter{getBasis}
\newcommand{\bsplinebasis}[4] % input: #1 knot vector, #2 order,  #3 id of basis function, #4 output
{						% example: \bsplinebasis{0&0&0&1&2&2&2}{2}{0}{\output}
    \newarray\vKnots 	% remark index start from 1
    \readarray{vKnots}{#1}

    \pgfmathparse{#3}  \let\i\pgfmathresult
    \pgfmathparse{#2}  \let\p\pgfmathresult

    \setcounter{countvalues}{0}
    \setcounter{getBasis}{\p} %\addtocounter{getBasis}{1}
    \pgfmathparse{int(\i+\p)}
    \foreach \knotspan in {\i,...,\pgfmathresult}
    {
	\pgfmathparse{ int(\knotspan+1+1) } \checkvKnots(\pgfmathresult)
	\pgfmathsetmacro{\tmpR}{\cachedata} 						%R: \tmpR
	
	\pgfmathparse{ int(\knotspan+1) } \checkvKnots(\pgfmathresult) 
	\pgfmathsetmacro{\tmpL}{\cachedata} 						%L: \tmpL

	\pgfmathparse{ \tmpR - \tmpL } \let\spansize\pgfmathresult  			%S: \spansize 
   
        \pgfmathparse{ \spansize > 0.0 } \let\bNonZero\pgfmathresult
        \ifthenelse{ \bNonZero = 1 }
        {
            \foreach \percentU in {0,10,...,100}
            {
		\pgfmathparse{\tmpL+\spansize*\percentU/100} \let\u\pgfmathresult

		\bsplinevalue{#1}{\p}{\u}{\knotspan}{\Basis} 
		\def\basisfuncarray{{\Basis}} 				% remark index start from 0
		
		\pgfmathparse{\basisfuncarray[\thegetBasis]} %\pgfmathresult
		
		\ifthenelse{\thecountvalues=0}
		{ 
			\xdef\nodeB{"\u,\pgfmathresult"}
		}{
			\xdef\nodeB{\nodeB,"\u,\pgfmathresult"}  
		}
		\addtocounter{countvalues}{1}
            }
        }{}
        \addtocounter{getBasis}{-1}
    }
    
	\xdef#4{\nodeB}

	\delarray\vKnots
}

\newcommand{\bsplinebasisspan}[6] % input: 	#1 knot vector, #2 order,  #3 segment id of basis function within (#4), 
{							% 		#4 id of knot span,  #5 id of plotted span, #6 output
    \newarray\vKnots 	% remark index start from 1
    \readarray{vKnots}{#1}

    \pgfmathparse{#5}  \let\plotknotspan\pgfmathresult
    \pgfmathparse{#4}  \let\splineknotspan\pgfmathresult
    \pgfmathparse{#3}  \let\i\pgfmathresult
    \pgfmathparse{#2}  \let\p\pgfmathresult

    \setcounter{countvalues}{0}
    \setcounter{getBasis}{\i} 
    \foreach \knotspan in {\plotknotspan}
    {
	\pgfmathparse{ int(\knotspan+1+1) } \checkvKnots(\pgfmathresult)
	\pgfmathsetmacro{\tmpR}{\cachedata} 						%R: \tmpR
	
	\pgfmathparse{ int(\knotspan+1) } \checkvKnots(\pgfmathresult) 
	\pgfmathsetmacro{\tmpL}{\cachedata} 						%L: \tmpL

	\pgfmathparse{ \tmpR - \tmpL } \let\spansize\pgfmathresult  			%S: \spansize 
   
        \pgfmathparse{ \spansize > 0.0 } \let\bNonZero\pgfmathresult
        \ifthenelse{ \bNonZero = 1 }
        {
            \foreach \percentU in {0,10,...,100}
            {
		\pgfmathparse{\tmpL+\spansize*\percentU/100} \let\u\pgfmathresult

		\bsplinevalue{#1}{\p}{\u}{\splineknotspan}{\Basis} 
		\def\basisfuncarray{{\Basis}} 				% remark index start from 0
		
		\pgfmathparse{\basisfuncarray[\thegetBasis]} %\pgfmathresult
		
		\ifthenelse{\thecountvalues=0}
		{ 
			\xdef\nodeB{"\u,\pgfmathresult"}
		}{
			\xdef\nodeB{\nodeB,"\u,\pgfmathresult"}  
		}
		\addtocounter{countvalues}{1}
            }
        }{}
        \addtocounter{getBasis}{-1}
    }
    
	\xdef#6{\nodeB}

	\delarray\vKnots
}

\newcommand{\plotbsplinebasis}[4] 	% input: #1 knot vector, #2 order,  #3 id of basis function, #4 plot settings
{							% example: \plotbsplinebasis{0&0&0&1&2&2&2}{2}{0}{green,smooth}
	\bsplinebasis{#1}{#2}{#3}{\nodeOut}
	\def\nodearray{{\nodeOut}}

	\xdef\name{ }
	\addtocounter{countvalues}{-1}
	\foreach \i in {0,...,\thecountvalues}
	{
		\pgfmathparse{\nodearray[\i]}
		\coordinate (point\i) at (\pgfmathresult);	  
		\xdef\name{ \name (point\i)  }
	}
	
	\draw[#4] plot coordinates{ \name };
	
	\xdef\name{ }
}

\newcommand{\plotbsplinesegment}[6] 	% input: 	#1 knot vector, #2 order,  #3 segment id of basis function within (#4),  
{								%		#4 id of knot span, #5 id of plotted span, #6 plot settings
	\bsplinebasisspan{#1}{#2}{#3}{#4}{#5}{\nodeOut}
	\def\nodearray{{\nodeOut}}

	\xdef\name{ }
	\addtocounter{countvalues}{-1}
	\foreach \i in {0,...,\thecountvalues}
	{
		\pgfmathparse{\nodearray[\i]}
		\coordinate (point\i) at (\pgfmathresult);	  
		\xdef\name{ \name (point\i)  }
	}
	
	\draw[#6] plot coordinates{ \name };
	
	\xdef\name{ }
}

\newcommand{\plotbsplineaccumulated}[5] 	% input: 	#1 knot vector, #2 order,  #3 subdivision coefficient,  
{						%		#4 id of knot span,  #5 plot settings
    \newarray\vKnots 	% remark index start from 1
    \readarray{vKnots}{#1}
    \newarray\vSubCoef 	% remark index start from 1
    \readarray{vSubCoef}{#3}
    
    \pgfmathparse{#4}  \let\plotknotspan\pgfmathresult
    \pgfmathparse{#4}  \let\splineknotspan\pgfmathresult
    \pgfmathparse{#2}  \let\p\pgfmathresult
    
    \setcounter{countvalues}{0}
    \pgfmathparse{int( \p+1) } \let\lastIndex\pgfmathresult
    \foreach \knotspan in {\plotknotspan}
    {
        \pgfmathparse{ int(\knotspan+1+1) } \checkvKnots(\pgfmathresult)
        \pgfmathsetmacro{\tmpR}{\cachedata} 						%R: \tmpR
        
        \pgfmathparse{ int(\knotspan+1) } \checkvKnots(\pgfmathresult) 
        \pgfmathsetmacro{\tmpL}{\cachedata} 						%L: \tmpL
        
        \pgfmathparse{ \tmpR - \tmpL } \let\spansize\pgfmathresult  			%S: \spansize 
        
        \pgfmathparse{ \spansize > 0.0 } \let\bNonZero\pgfmathresult
        \ifthenelse{ \bNonZero = 1 }
        {
            \foreach \percentU in {0,10,...,100}
            {
                \pgfmathparse{\tmpL+\spansize*\percentU/100} \let\u\pgfmathresult
                
                \bsplinevalue{#1}{\p}{\u}{\splineknotspan}{\Basis} 
                \def\basisfuncarray{{\Basis}} 				% remark index start from 0
                
                \setcounter{getBasis}{0} 
                \pgfmathparse{\basisfuncarray[\thegetBasis]} %\pgfmathresult
                \let\basisValue\pgfmathresult
                
                \checkvSubCoef(1) \pgfmathsetmacro{\coef}{\cachedata}  
                \pgfmathparse{ \basisValue * \coef } \let\first\pgfmathresult
                
                \xdef\lastx{\first}
                \foreach \i in {2,...,\lastIndex}
                { 
                    \addtocounter{getBasis}{1}       
                    \pgfmathparse{\basisfuncarray[\thegetBasis]}
                    \let\basisValue\pgfmathresult
                    
                    \checkvSubCoef(\i) \pgfmathsetmacro{\coef}{\cachedata}  
                    \pgfmathparse{ \lastx + \basisValue * \coef } \let\sum\pgfmathresult
                    
                    \xdef\lastx{\sum}
                    
                    \ifthenelse{\i=\lastIndex}
                    {
                        \ifthenelse{\thecountvalues=0}
                        { 
                            \xdef\nodeBB{"\u,\sum"}
                        }{
                            \xdef\nodeBB{\nodeBB,"\u,\sum"}  
                        }
                        \addtocounter{countvalues}{1}
                    }{}
                    
                }
            }
        }{}
    }
    
    \delarray\vKnots
    \delarray\vSubCoef
    
    \def\nodearray{{\nodeBB}}
    
    \xdef\name{ }
    \addtocounter{countvalues}{-1}
    \foreach \i in {0,...,\thecountvalues}
    {
        \pgfmathparse{\nodearray[\i]}
        \coordinate (point\i) at (\pgfmathresult);                
        \xdef\name{ \name (point\i)  }
    }
    
    \draw[#5] plot coordinates{ \name };
    
    \xdef\name{ }
}

\begin{document}
    
\title{Efficient simulation of inclusions and reinforcement bars with the isogeometric Boundary Element method}
\begin{frontmatter}

\author[ifbaddr]{Gernot Beer\corref{cor1}}
\author[ifbaddr]{Christian Dünser}
\author[sun]{Eugenio Ruocco}
\author[unife]{Vincenzo Mallardo}

\address[ifbaddr]{Institute of Structural Analysis, Graz University
  of Technology, Lessingstraße 25/II, 8010 Graz, Austria}

\address[sun]{Department of Engineering, University of Campania "L. Vanvitelli", Via Roma 28, 81031 Aversa, Caserta,
Italy}

\address[unife]{Department of Architecture, University of Ferrara, Via Quartieri 8, 44121 Ferrara, Italy}

\cortext[cor1]{Corresponding author.
  Tel.: +43 316 873 6181, fax: +43 316 873 6185, mail: \url{gernot.beer@tugraz.at}, web: \url{www.ifb.tugraz.at}}

\begin{abstract}
The paper is concerned with the development of efficient and accurate solution procedures for the isogeometric boundary element method (BEM) when applied to problems that contain inclusions that have elastic properties different to the computed domain.
This topic has been addressed in previous papers but the approach presented here is a considerable improvement in terms of efficiency and accuracy. 

One innovation is that initial stresses instead of body forces are used. This then allows a one step solution without iteration.
In addition, a novel approach is used for the computation of strains, that avoids the use of highly singular fundamental solutions.
Finally, a new type of inclusion is presented that can be used to model reinforcement bars or rock bolts and where analytical integration can be used. Test examples, where results are compared with Finite Element simulations, show that the proposed approach is sound.

\end{abstract}

\end{frontmatter}

\section{Introduction}
Since the publication of the first paper on the topic \cite{Hughes2005a}, isogeometric analysis has gained increased popularity. 
The majority of applications have been with the Finite Element method  (FEM) and much less with the Boundary Element method.
However, the advantage of the BEM, that only the discretisation of the boundary is required, makes it an ideal companion to CAD.
First applications of the isogeometric BEM were published in elasticity in 2-D \cite{simpson2012two,simpson2013isogeometric} and in 3-D \cite{scott2013isogeometric}. The method was also applied to problems in acoustics \cite{simpson2014acoustic} and electromagnetic scattering \cite{simpson2018isogeometric}. Ways of accelerating isogeometric BEM solutions were published \cite{simpson2016acceleration}.
In \cite{Marussig2015} the concept of a geometry independent field approximation, which involved a decoupling of the geometry definition and the approximation of the unknown,  was first introduced. This concept that was later adopted by others \cite{doi:10.1002/nme.5778}. The seamless integration of BEM and CAD was discussed in \cite{marussig2016b}.
The method was applied to viscous flow in \cite{Beer2018,Beer17a}.
In a recent book published on the isogeometric BEM \cite{BeerMarussig} it was shown how geometrical information can be taken directly from CAD data and that efficient and accurate simulations with very few unknowns can be obtained.

One fact that has hampered the widespread use of the isogeometric BEM is that fundamental solutions, on which the method is based, exist only for elastic material properties and homogeneous domains. However, as shown in the cited literature, this obstacle can be overcome by the use of body forces.
The concept of body forces has already been applied to the solution of problems with elasto-plastic inclusions in \cite{Beer17,Beer2016}. However, the approach presented there involved an iterative process and the conversion of initial stresses to body forces, resulting in a lack of efficiency.  Here we address this problem and propose remedies. In addition a new type of inclusion, namely reinforcement bars or rock bolts, is introduced. For this type of inclusion analytical integration can be used, resulting in a significant increase in efficiency. 

\section{Theory}
As in previous work \cite{Beer17} we propose to solve the problem by considering initial stresses generated inside the inclusions.
The initial stresses $ \myVecGreek{\sigma}_{0}^{e} $ due to a difference in elastic properties are given by:
\begin{equation}
\label{InitialS}
\myVecGreek{\sigma}_{0}^{e} =( \mathbf{ D} - \mathbf{ D}_{incl}) \myVecGreek{\epsilon} 
\end{equation}
where $\myVecGreek{\epsilon} $ is the total strain, $ \mathbf{ D}$ and $ \mathbf{ D}_{incl}$ is the elasticity matrix for the domain and the inclusion, respectively. 

We may also consider initial stresses $ \myVecGreek{\sigma}_{0}^{p} $ due to elasto-plasticity, given by:
\begin{equation}
\label{InitialSp}
\myVecGreek{\sigma}_{0}^{p} =( \mathbf{ D} - \mathbf{ D}_{ep, incl}) \myVecGreek{\epsilon} 
\end{equation}
where  $ \mathbf{ D}_{ep,incl}$ is the elasto-plastic matrix of the inclusion. 

In \cite{Beer17}  the initial stresses were converted to body forces and the solution was obtained by iteration. We now propose to use the initial stresses directly and to solve the elastic problem without iterations. There are several advantages to this approach, as will be discussed, but one benefit is that it allows the modelling of non-linear material behaviour for the case where the inclusion has different elastic properties to the domain. The previously published iterative method was not suitable to be combined with non-linear iterations.
Another novel part is that (in addition to the mapping methods published in \cite{BeerMarussig}) we develop special mapping methods for reinforcement bars and rock bolts.
In the following we first establish the governing integral equations and then discuss in detail how the arising volume and surface integrals are evaluated.

\subsection{Governing integral equations}
Consider a domain $\domain$ with a boundary $\boundary$, containing a subdomain $\domain_{0}$ where initial stresses $ \myVecGreek{\sigma}_{0}(\fieldpt)$ are present.
We apply the theorem by Betti and the collocation method to arrive at the governing integral equations.
This means that we set the work done on the boundary $\boundary$ by tractions $\fund{T}$ times displacements $\myVec{u}$ equal to the work done by displacements $\fund{U}$ times tractions $\myVec{t}$. We assume $\fund{T}$ and $\fund{U}$ to be fundamental solutions of the governing differential equation at $\fieldpt$ due to a source at $\sourcept_{n}$ and $\myVec{u}$, $\myVec{t}$ to be boundary values.
If initial stresses are present, additional work is done in the domain $\domain_{0}$ by the initial stresses $ \myVecGreek{\sigma}_{0}(\fieldpt)$ times strains $\fund{E} (\sourcept_{n},\fieldpt)$.
Therefore the regularised integral equation is written as :
\begin{equation}
\label{eq_displacement}
    \begin{aligned}
% \begin{eqnarray}
%\label{BEM}
\int_{\boundary} \fund{T}(\sourcept_{n},\fieldpt) ( \myVec{u}(\fieldpt) - \myVec{u}(\sourcept_{n})) \ d\boundary(\fieldpt) - \mathbf{ A}_{n} \myVec{u}(\sourcept_{n}) &=& \int_{\boundary} \fund{U}(\sourcept_{n},\fieldpt) \ \myVec{t}(\fieldpt) \ d\boundary(\fieldpt) \\
% \nonumber
+  \int_{\domain_{0}} \fund{E} (\sourcept_{n},\fieldpt)  \myVecGreek{\sigma}_{0} (\fieldpt) d \domain_{0} (\fieldpt)  .
% \end{eqnarray}
\end{aligned}
\end{equation}
where $\mathbf{ A}_{n}$ represents the azimuthal integral that is zero for finite domain problems. The derivation of Eq. (\ref{eq_displacement}) and the fundamental solutions $\fund{U}$ und $\fund{T}$ are presented in \cite{BeerMarussig}.

The fundamental solution $\fund{E}$ is given by:
\begin{equation}
\label{ }
E_{ijk}= \frac{-C}{r^{2}}\left[C_{3}(r_{,k}\delta_{ij} + r_{,j}\delta_{ik}) - r_{,i}\delta_{jk} + C_{4} \ r_{,i} r_{,j} r_{,k}\right]
\end{equation}
where $r=|\fieldpt - \sourcept|$, $r_{,i}=  \frac{r_{i}}{r}$ and $\delta_{ij}$ is the \textit{Kronecker Delta}. The constants are: $C=\frac{1}{16 \pi G (1- \nu})$, $C_{3}=1 - 2\nu$ and $C_{4}=3$ where $G$ is the shear modulus and $\nu$ the Poisson's ratio.

In the following we will use matrix algebra and it is therefore necessary to convert the stress and strain tensors to matrices, using \textit{Voigt} notation.

The initial stresses can be written as:
\begin{equation}
\label{Voightnot}
\myVecGreek{\sigma}_{0}= \left \{\begin{array}{c}\sigma_{11} \\\sigma_{22} \\\sigma_{33} \\\sigma_{12} \\\sigma_{23}\\\sigma_{13}\end{array}\right\}_{0}
\end{equation}
The tensor $E_{ijk}$ is converted to a matrix $\fund{E}$:
\begin{equation}
\label{eq_Voigt_convertion}
\fund{E}= \left[\begin{array}{cccccc}E_{111} & E_{122} & E_{133}  & E_{112}  + E_{121}& E_{123} + E_{132} & E_{113} +E_{131}  \\E_{211} & E_{222} & E_{233}  & E_{212} + E_{221}  & E_{223}  + E_{232} & E_{213} + E_{231}\\E_{311} & E_{322} & E_{333}  & E_{312} + E_{321} & E_{323}  + E_{332}& E_{313} + E_{331}\end{array}\right]
\end{equation}

\subsection{Discretised integral equations}
\subsubsection{Geometry approximation of the boundary}
For the geometrical discretisation we divide the boundary of the problem into patches and the inclusion into subdomains.
After the discretisation into patches we obtain the following boundary element equations (omitting the volume term for the moment):
\begin{equation}
    \begin{aligned}
        \label{CollocR1Second}
        \sum_{e=1}^{E} \int_{\Gamma_{e}}  \fund{U}(\sourcept_{n},\fieldpt) \ \myVec{\dual}^{e}(\fieldpt)\ d\Gamma_{e}(\fieldpt)  
        &- \sum_{e=1}^{E} \int_{\Gamma_{e}} \fund{T}(\sourcept_{n},\fieldpt)  \myVec{\primary}^{e}(\fieldpt) d \Gamma_{e}  \\
        &+  \left[\sum_{e=1}^{E} \ \left(\int_{\Gamma_{e}} \fund{T}(\sourcept_{n},\fieldpt) d \Gamma_{e} \right)  - \mathbf{A}_{n}\right] \myVec{\primary}(\sourcept_{n})= 0
    \end{aligned}
\end{equation}
where $e$ specifies the patch number and $E$ is the total number of patches.
The geometry of patches is defined using NURBS basis functions $\NURBS_{i}$ of local coordinates $\xi , \eta$:
\begin{equation}
%\label{ }
\pt{x} (\xi,\eta)= \sum_{i=1}^{I} \NURBS_{i} (\xi, \eta) \pt{x}_{i}.
\end{equation}
where $I$ is the total number of control points and $\pt{x}_{i}$ are control point coordinates.

\subsubsection{Approximation of boundary values}
For the approximation of the boundary values we use a geometry independent field approximation.
Unknown values are approximated by
\begin{equation}
\begin{aligned}
\label{eq6:ApproxSecond}
 \myVec{\hat{\primary}}^{e} (\xi,\eta) &=  \sum_{k=1}^{K}   \hat{ \NURBS}_{k} (\xi,\eta) \   \myVec{\hat{\primary}}_{k}^{e}\\
 \myVec{\hat{\dual}}^{e}(\xi,\eta)  &=   \sum_{k=1}^{K}    \hat{\NURBS}_{k} (\xi,\eta) \  \myVec{\hat{\dual}}_{k}^{e}.
\end{aligned}
\end{equation}
and may be refined using the standard methods used in isogeometric analysis such as knot insertion and order elevation.
Our refinement philosophy is to take the NURBS functions that define the geometry of the problem and refine them as necessary.

Known values are defined by
\begin{equation}
\begin{aligned}
\label{}
 \myVec{\bar{\primary}}^{e}(\xi,\eta)  &=  \sum_{k=1}^{\bar{K}}    \bar{\NURBS}_{k}(\xi,\eta) \   \myVec{\bar{\primary}}_{k}^{e}\\
 \myVec{\bar{\dual}}^{e}(\xi,\eta)  &=   \sum_{k=1}^{\bar{K}}    \bar{\NURBS}_{k}(\xi,\eta) \  \myVec{\bar{\dual}}_{k}^{e}.
\end{aligned}
\end{equation}
where $\hat{\NURBS}_{k}$ and $\bar{\NURBS}_{k}$ are basis functions, which may be different from the ones defining the geometry and the unknown values. $\myVec{\hat{\primary}}_{k}^{e}, \myVec{\hat{\dual}}_{k}^{e}$ are the parameter values at control points.\footnote{It should be noted here that for known values that have a simple variation we may even use lower order Lagrange polynomials for $\bar{\NURBS}_{k}$, in which case $\myVec{\bar{\primary}}_{k}^{e}$ or $ \myVec{\bar{\dual}}_{k}^{e}$ represent nodal values.} 

\subsection{Geometry definition of inclusions}

For the simulation we use two different types of inclusion: General inclusion to represent a volume and a linear inclusion to represent reinforcement bars.

\subsubsection{General inclusion}
\begin{figure}[h]
\begin{center}
\begin{overpic}[scale=0.7]{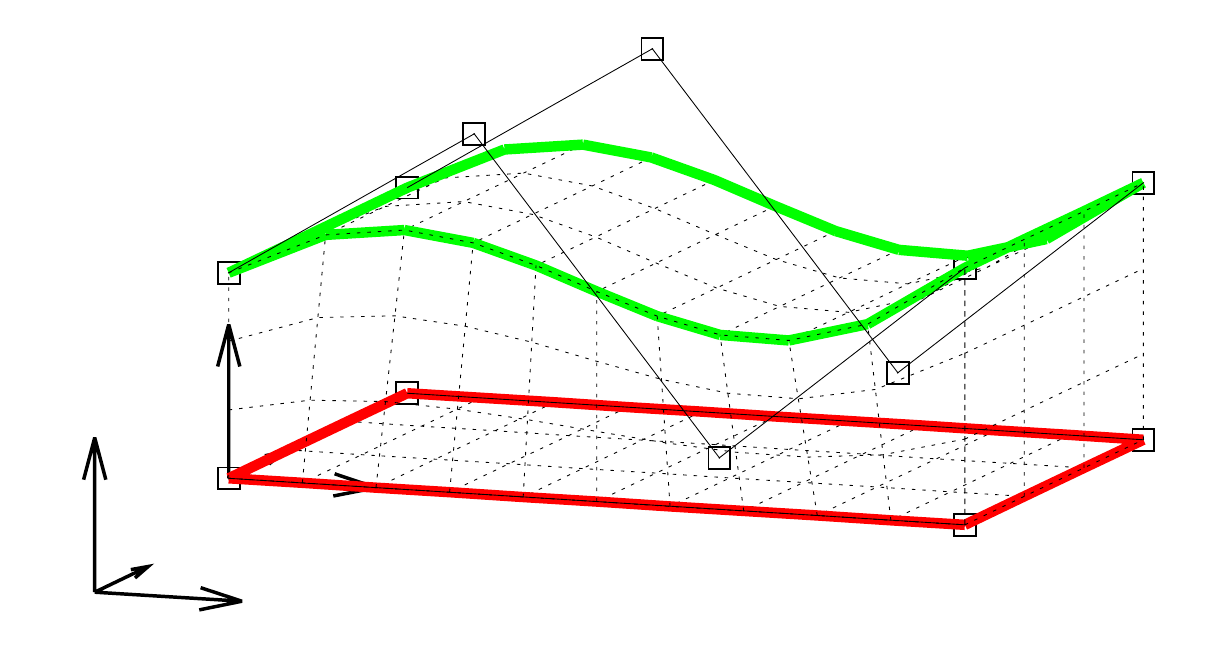}
 \put(20,5){$x$}
 \put(6,20){$z$}
 \put(10,10){$y$}
  \put(20,28){$r$}
 \put(30,15){$s$}
 \put(23,18){$t$}
\end{overpic}
\begin{overpic}[scale=0.5]{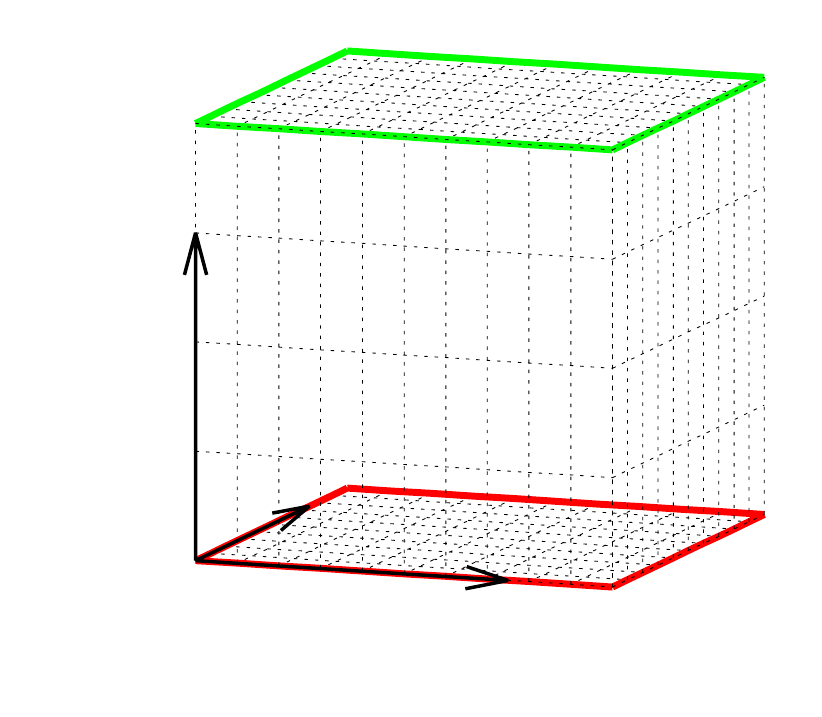}
  \put(60,20){$s$}
 \put(35,30){$t$}
 \put(15,60){$r$}
\end{overpic}
\caption{Mapping of 3-D inclusion showing the bottom and top NURBS surfaces and the associated control points defining the inclusion: Left in global $\pt{x}$, right in local $\pt{s}$ space. Also shown are subregions for the volume integration.}
\label{fig:Incl3D}
\end{center}
\end{figure}
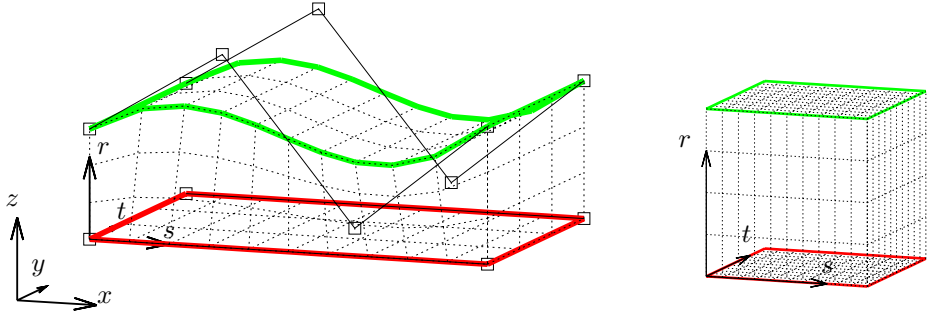

General inclusions are defined by bounding NURBS surfaces.
We establish a local coordinate system $\pt{s}=(s,t,r)^{\trans}=[0,1]^3$ as shown in \myfigref{fig:Incl3D} and map from local $\pt{s}$ coordinates to global $\pt{x}$ coordinates.
The global coordinates of a point $\pt{x}$ with the local coordinates $\pt{s}$ are given by
\begin{equation}
%\label{ }
\pt{ x}({s,t,r})= (1-r) \ \pt{ x}^{I}(s,t) + {r} \ \pt{ x}^{II}({s,t})
\end{equation}
where
\begin{align}
%\label{ }
  \pt{ x}^{I}({s,t})=\sum_{k=1}^{K^{I}} R_{k}^{I}({s,t}) \ \pt{ x}_{k}^{I} && \text{and} &&\pt{ x}^{II}({s,t})=\sum_{k=1}^{K^{II}} R_{k}^{II}({s,t}) \ \pt{ x}_{k}^{II} .
\end{align}
The superscript $I$ relates to the bottom (red) surface and $II$ to the top (green) bounding surface and $ \pt{ x}_{k}^{I} $, $ \pt{ x}_{k}^{II} $ are control point coordinates. $K^{I}$ and $K^{II}$ represent the number of control points, $R_{k}^{I}({s,t})$ and $R_{k}^{II}({s,t})$ are NURBS basis functions. Note that there is a one to one mapping between the local surface coordinates $\uu,\vv$ and the local coordinates $s,t$.

The derivatives are given by
\begin{equation}
%\label{ }
\begin{aligned}
  \frac{\partial \pt{ x}({s,t,r})}{\partial {s}}&=& (1-{r}) \ \frac{\partial \pt{ x}^{I}({s,t})}{\partial {s}} &+&& {r} \ \frac{\partial \pt{ x}^{II}({s,t}) }{\partial {s}} \\
    \frac{\partial \pt{ x}({s,t,r})}{\partial {t}}&=& (1-{r}) \ \frac{\partial \pt{ x}^{I}({s,t})}{\partial {t}} &+&& {r} \ \frac{\partial \pt{ x}^{II}({s,t}) }{\partial {t}} \\
  \frac{\partial \pt{ x}({s,t,r})}{\partial {r}}&=& -\pt{ x}^{I}({s,t}) &+&&
  \ \pt{ x}^{II}({s,t})
\end{aligned}
\end{equation}
where for example:
\begin{align}
  \frac{\partial\pt{ x}^{I}({s,t})}{\partial {s}}=\sum_{k=1}^{K^{I}} \frac{\partial R_{k}^{I}({s,t})}{\partial {s}} \ \pt{ x}_{k}^{I} &&\text{and} &&
  \frac{\partial\pt{ x}^{II}({s,t})}{\partial {}s}=\sum_{k=1}^{K^{II}}
  \frac{\partial R_{k}^{II}({s,t})}{\partial {s}} \ \pt{ x}_{k}^{II}  .
\end{align}
The Jacobi matrix of this mapping is
\begin{equation}
\label{eq:jac3D}
\mathbf{J}=
\begin{pmatrix}
  \frac{\partial \pt{x}}{\partial {s}}   \\ \\
  \frac{\partial \pt{x}}{\partial {t}}  \\ \\
  \frac{\partial \pt{x}}{\partial {r}}  
\end{pmatrix}
\end{equation}
and the Jacobian is $J=| \mathbf{ J} |$.

\subsubsection{Linear inclusion, reinforcement bar/rock bolt}
This type of inclusions is used for reinforcement bars and rock bolts. Here we assume that the geometry is defined by a linear NURBS curve and that the bar has a circular cross-section with radius $R$ over which the stress and strain are assumed constant. \emph{The assumption is that the area of the cross-section of the inclusion is significantly smaller that of the medium it is embedded in, allowing simplifications to be introduced for the integration.}

\begin{remark}
It should be noted here that a higher order NURBS could also be used to define the geometry of a curved bar. The restriction to linear NURBS has been imposed in this paper because it allows the integration to be carried out analytically, resulting in an efficient simulation if many rock bolts are present.
\end{remark}

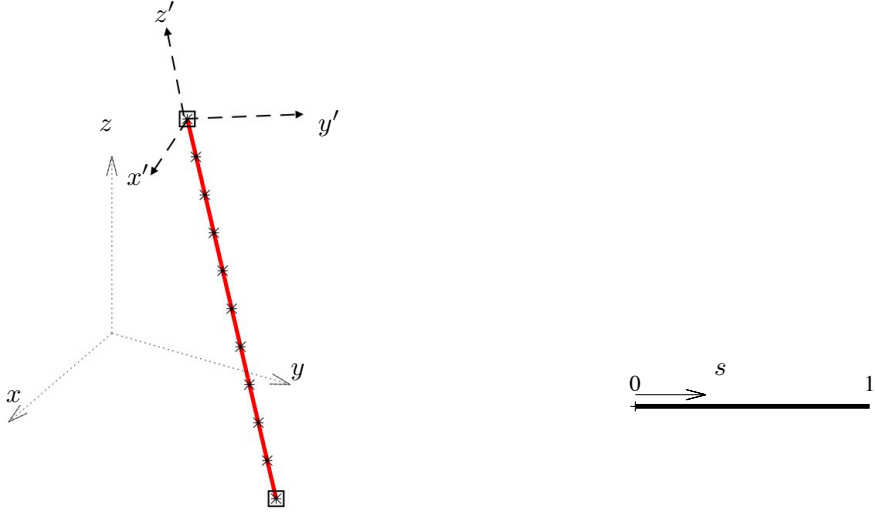
\begin{figure}
\begin{center}
\begin{overpic}[scale=0.6]{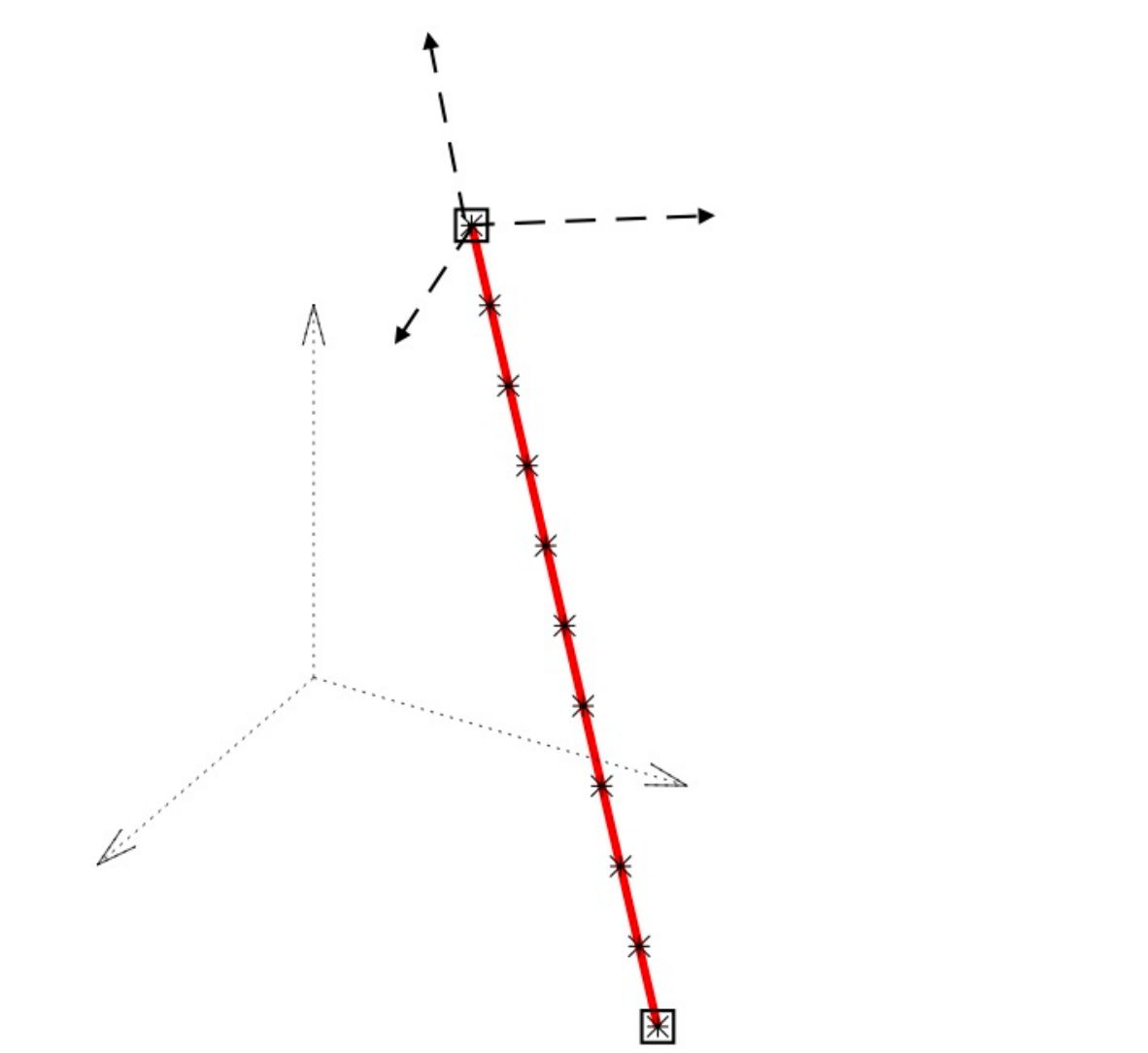}
\put(25,70){$z$}
\put(8,20){$x$}
\put(60,25){$y$}
\put(65,70){$y^{\prime}$}
\put(30,60){$x^{\prime}$}
\put(35,90){$z^{\prime}$}
\end{overpic}
\begin{overpic}[scale=0.6]{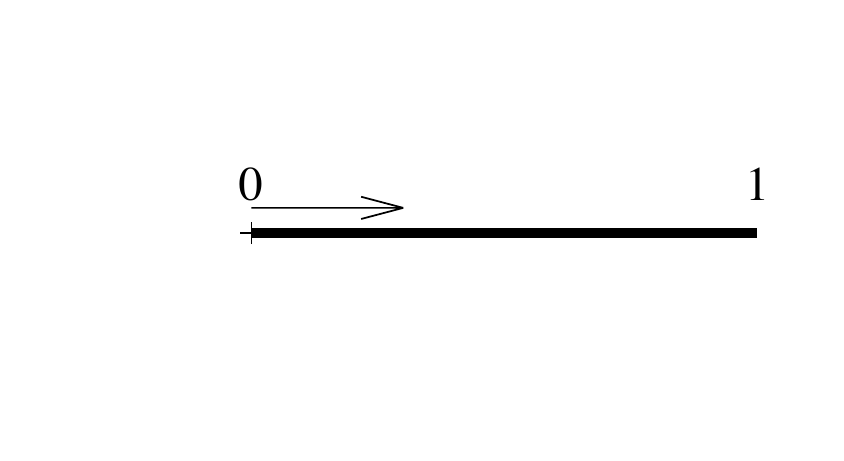}
\put(50,35){$s$}
\end{overpic}
\caption{Definition of linear inclusion by a NURBS curve with control points as hollow squares in global (left) and local (right) coordinates. 
Local axes $x^{\prime}, y^{\prime}, z^{\prime}$ are shown at the end of the bar. Also shown are internal points.}
\label{Rebar}
\end{center}
\end{figure}

We establish a local coordinate system $s=\left[0,1\right]$ as shown on the right in \myfigref{Rebar}.
The global coordinates of a point $\pt{x}$ with the local coordinate $s$ are given by
\begin{eqnarray}
%\label{ }
  \pt{ x} ({s})=\sum_{k=1}^{K} R_{k}(s) \ \pt{ x}_{k}  
\end{eqnarray}
where $K$ is the number of control points, $R_{k}(s)$ are NURBS basis functions and $\pt{ x}_{k}$ are control point coordinates.
We also define a local coordinate system whereby the $z^{\prime}$ axis is along the bar specified by unit vector $\mathbf{ v}_{z^{\prime}}$. 
The vector along the bar is given by
\begin{eqnarray}
\mathbf{ V}_{z^{\prime}}= \frac{\partial\pt{ x}({s})}{\partial {s}}&=&\sum_{k=1}^{K} \frac{\partial R_{k}({s})}{\partial {s}} \ \pt{ x}_{k} 
\end{eqnarray}
The Jacobian is
\begin{equation}
\label{ }
J= \sqrt{\mathrm{V}_{z^{\prime}_x}^{2} + \mathrm{V}_{z^{\prime}_y}^{2} + \mathrm{V}_{z^{\prime}_z}^{2}}
\end{equation}
The unit vector in $z^{\prime}$ direction is given by
\begin{equation}
\label{ }
\mathbf{ v}_{z^{\prime}}= \mathbf{ V}_{z^{\prime}}/J
\end{equation}
The direction of the other axes are constructed perpendicular to $z^{\prime}$, as will be shown later.

\subsection{Evaluation of integrals}
The surface integrals in Eq. (\ref{CollocR1Second}) are evaluated numerically using established procedures as outlined in  \cite{BeerMarussig}. 
The evaluation of the volume integrals will be discussed later.

\subsection{System of equations}

The discretised integral equations can be written as:
\begin{equation}
\label{DisIE}
[\mathbf{ L}] \{\mathbf{ x}\} = \{\mathbf{ r}\} + [\mathbf{ B}_{0}] \{ \myVecGreek{\sigma}_{0}\} 
\end{equation}
where  $\mathbf{ L}$ is the left hand side, $\mathbf{ x}$ is the vector of unknowns and $\mathbf{ r}$ is the right hand side (for details of derivation see\cite{BeerMarussig}), $\{ \myVecGreek{\sigma}_{0}\}$ is a vector that gathers all initial stress components at grid points inside the inclusions.

For the numerical integration, used for the general inclusions, we need the values of the initial stress at Gauss points. The value of initial stress at a point with the local coordinates $\pt{s}=(s,t,r)^\mathrm{T}=[0,1]^3$ is obtained by interpolation between grid points.
\begin{equation}
\label{eq11:Interpol}
 \myVecGreek{\sigma}_{0}(\pt{s})= \sum_{j=1}^{J} M_{j}^{\sigma}(\pt{s})  \myVecGreek{\sigma}_{0j}
\end{equation}
where $  \myVecGreek{\sigma}_{0j} $ is the initial stress vector at grid point $j$ with the local coordinate $\pt{s}_{j}$. $M_{j}^{\sigma}(\pt{s})$ are linear or constant basis functions, which will be shown later. 

The sub-matrices of matrix $ [\mathbf{ B}_{0}]$, related to collocation point $n$ and grid point $j$, are given by
\begin{equation}
\label{ }
\mathbf{ B}_{0nj} = \int_{\domain_{0}} \fund{E} (\sourcept_{n},\fieldpt) M_{j}^{\sigma} (\fieldpt) d \domain_{0} (\fieldpt) 
\end{equation}

\subsection{Computation of values at grid points inside the inclusion}
\subsubsection{Computation of displacements}
For the solution we need to compute the displacements and strain at points inside the inclusion.
The displacement $ \mathbf{ u} $ at a point $\pt{x}$ inside the inclusion is given by:
\begin{equation}
    \begin{aligned}
        \label{displ}
       \mathbf{ u}(\pt{x}) &= \int_{\boundary} \left[ \fund{U}(\pt{x},\fieldpt) \   \myVec{\dual} (\fieldpt)  - \fund{T}(\pt{x},\fieldpt) \   \myVec{\primary} (\fieldpt) \right] d\boundary (\fieldpt) \\
%         \nonumber
        &+ \int_{\domain_{0}} \fund{E} (\pt{x},\fieldpt) \myVecGreek{\sigma}_{0} (\fieldpt)  d \domain_{0} (\fieldpt)     
    \end{aligned}
\end{equation}

We gather displacement vectors at all grid points in a vector $\{\mathbf{ u} \}$. and obtain:
\begin{equation}
\label{ }
\{\myVec{u}\}= [\hat{\mathbf{ A}}] \mathbf{ x} + \{\bar{\mathbf{ c}}\} + [\bar{\mathbf{ B}}_{0}]\{ \myVecGreek{\sigma}_{0}\} 
\end{equation}
where $[\hat{\mathbf{ A}}]$ is an assembled matrix that multiplies with the unknown $\mathbf{ x}$ and $\{\bar{\mathbf{ c}}\}$ collects the displacement contribution due to given BC's. $ [\bar{\mathbf{ B}}_{0}]$ is similar to $ [\mathbf{ B}_{0}]$ except that the grid point coordinates $\pt{x}_{i}$ replace the source point coordinates $\sourcept_{n}$.

Because of the singularity of $\fund{T}$ the displacements can not be computed on the problem boundary. So if the inclusion intersects a boundary patch we recover the displacement from the boundary values on the patch.
For points on a patch boundary ($\pt{x}_{k}$) we replace  Eq. (\ref{displ}) by:
\begin{equation}
\label{ }
\mathbf{ u}(\pt{x}_{k})=\sum_{i}^{I} R_{i}^{u}(\xi_{k},\eta_{k}) \mathbf{ u}_{i}^{e}
\end{equation}
where $R_{i}^{u}(\xi,\eta)$ are the NURBS basis functions used for approximating the displacements in patch $e$, that contains the point $\pt{x}_{k}$ and $\xi_{k},\eta_{k}$ are the local coordinates of the point.
The matrix $[\hat{\mathbf{ A}}]$ and the vector $\{\bar{\mathbf{ c}}\} $ have to be modified for these grid points, whereas $ [\bar{\mathbf{ B}}_{0}]$ will contain zero rows.

The strain at grid points can be computed using derived fundamental solutions but because of their high singularity this would involve complicated integration schemes (see for example \cite{Gao2011}). We choose a simpler alternative borrowed from the FEM community by  using strain recovery, i.e. by taking the derivative of the displacements. 

We interpolate the displacements between grid points and obtain for the displacement at a point with the local coordinate $\pt{s}$:
\begin{equation}
\label{Interpolu}
 \myVec{u}(\pt{s})= \sum_{j=1}^{J} M_{j}(\pt{s})  \myVec{u}_{j}
\end{equation}
where $\myVec{u}_{j}$ is the displacement vector at grid point $j$ and $J$ is the number of grid points.  $M_{j}(\pt{s})$ are piecewise linear or parabolic  shape functions that will be shown later.

\subsubsection{Computation of strains, general inclusions}

 For general inclusions the strains are given by:
\begin{eqnarray}
\epsilon_{x} & = & \frac{\partial u_{x}}{\partial x} = \sum_{j=1}^{J} \frac{\partial M_{j}}{\partial x} u_{xj}\\
\epsilon_{y} & = & \frac{\partial u_{y}}{\partial y} = \sum_{j=1}^{J} \frac{\partial M_{j}}{\partial y} u_{yj}\\
\epsilon_{z} & = & \frac{\partial u_{z}}{\partial z} = \sum_{j=1}^{J} \frac{\partial M_{j}}{\partial z} u_{zj}\\
\gamma_{xy} & = & \frac{\partial u_{x}}{\partial y} + \frac{\partial u_{y}}{\partial x}= \sum_{j=1}^{J} \frac{\partial M_{j}}{\partial x} u_{yj} + \sum_{j=1}^{J} \frac{\partial M_{j}}{\partial y} u_{xj}\\
\gamma_{zy} & = & \frac{\partial u_{z}}{\partial y} + \frac{\partial u_{y}}{\partial z}= \sum_{j=1}^{J} \frac{\partial M_{j}}{\partial z} u_{yj} + \sum_{j=1}^{J} \frac{\partial M_{j}}{\partial y} u_{zj}\\
\gamma_{xz} & = & \frac{\partial u_{x}}{\partial z} + \frac{\partial u_{z}}{\partial x}= \sum_{j=1}^{J} \frac{\partial M_{j}}{\partial x} u_{zj} + \sum_{j=1}^{J} \frac{\partial M_{j}}{\partial z} u_{xj}
\end{eqnarray}
The strains gathered at grid point $k$ can be written in matrix notation:
\begin{equation}
\label{Bhat}
\myVecGreek{\epsilon}(\pt{x}_{k})= \hat{\mathbf{ B}}(\pt{x}_{k}) \{\myVec{u}\}
\end{equation}
where
\begin{equation}
\label{ }
\hat{\mathbf{ B}}(\pt{x}_{k})= \left(\begin{array}{ccc}\mathbf{ B}_{1} & \mathbf{ B}_{2} & \cdots\end{array}\right)
\end{equation}
and
\begin{equation}
\label{Bsub}
 \mathbf{ B}_{i}=\left(\begin{array}{ccc}\frac{\partial M_{i}}{\partial x} & 0 & 0 \\0 & \frac{\partial M_{i}}{\partial y} & 0 \\0 & 0 & \frac{\partial M_{i}}{\partial z}  \\ \frac{\partial M_{i}}{\partial y} & \frac{\partial M_{i}}{\partial x} & 0 \\ 0 & \frac{\partial M_{i}}{\partial z} & \frac{\partial M_{i}}{\partial y}  \\ \frac{\partial M_{i}}{\partial z} & 0 &\frac{\partial M_{i}}{\partial x} \end{array}\right)
\end{equation}

The global derivatives of $M_{j}$ are given by:
\begin{equation}
\label{ }
\left(\begin{array}{c}\frac{\partial M_{j}}{\partial x} \\ \\\frac{\partial M_{j}}{\partial y} \\ \\\frac{\partial M_{j}}{\partial z}\end{array}\right)= \mathbf{ J}^{-1} \left(\begin{array}{c}\frac{\partial M_{j}}{\partial s} \\ \\\frac{\partial M_{j}}{\partial t} \\ \\ \frac{\partial M_{j}}{\partial r}\end{array}\right)
\end{equation}
where $\mathbf{ J}$ is the Jacobian matrix Eq. (\ref{eq:jac3D}).  
For a linear inclusion we compute the strain in local directions as is shown later.

\paragraph{Definition of $M_{j}$}

For the interpolation between grid points we consider an equally spaced grid  in the local coordinate directions $s,t,r$.  The inclusion is divided then into regions of equal size ($\triangle s, \triangle t, \triangle r$) (if there are only 2 grid points in a direction only one region and linear interpolation is used). 
In the following we explain the interpolation in one direction ($s$) as the scheme is identical for the other directions ($t,r$).

\begin{figure}
\begin{center}
\begin{overpic}[scale=0.6]{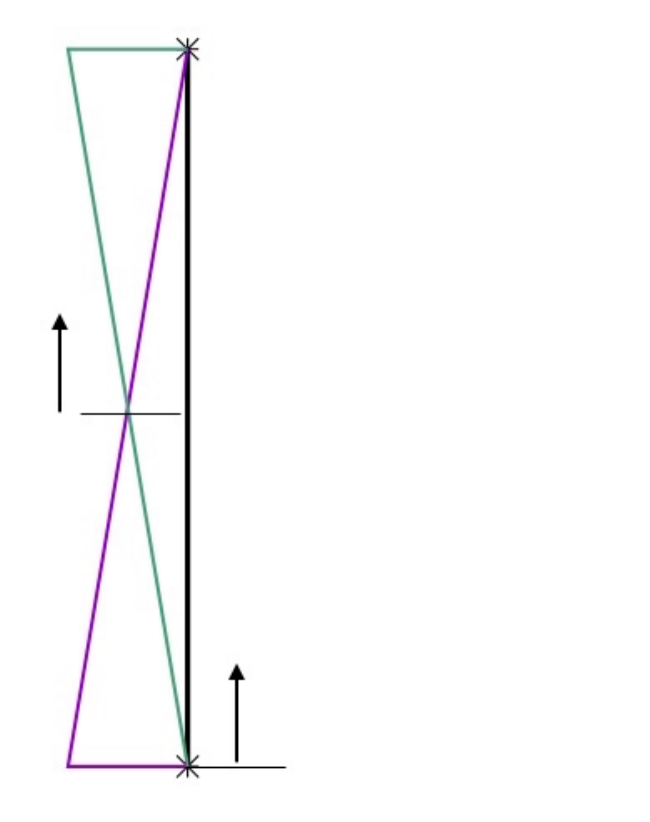}
\put(30,20){$s$}
\put(8,65){$\xi$}
\put(30,5){$j=1$}
\put(27,90){$j=2$}
\put(27,60){$\triangle s$}
\end{overpic}
\begin{overpic}[scale=0.6]{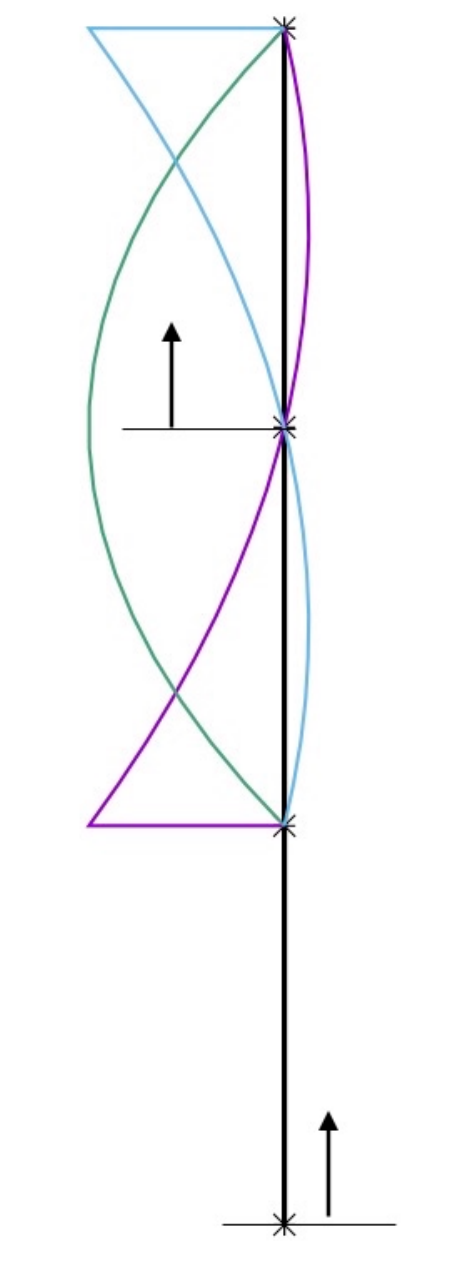}
\put(28,15){$s$}
\put(10,70){$\xi$}
\put(16,30){$s_{j}$}
\put(27,5){$j=1$}
\put(25,35){$j=2$}
\put(25,70){$j=3$}
\put(25,95){$j=4$}
\put(27,60){$2\triangle s$}
\end{overpic}
\begin{overpic}[scale=0.9]{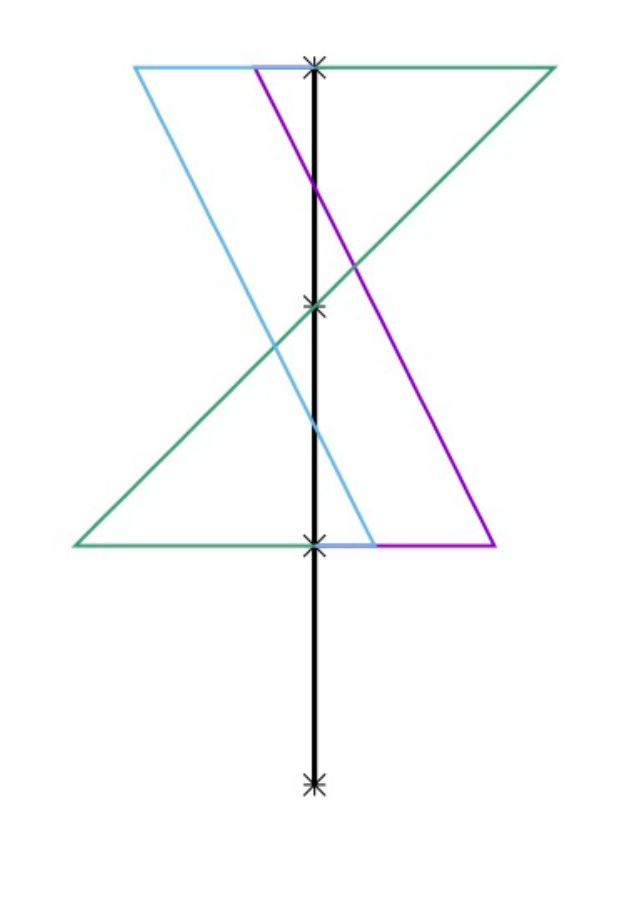}
\end{overpic}
\caption{Interpolation functions $\bar{M}_{j}(\xi)$ in $s$ direction for the case where there are only 2 internal points (left) and more than 2 (right). Also shown are the derivatives of the quadratic functions.}
\label{Parainterpol}
\end{center}
\end{figure}
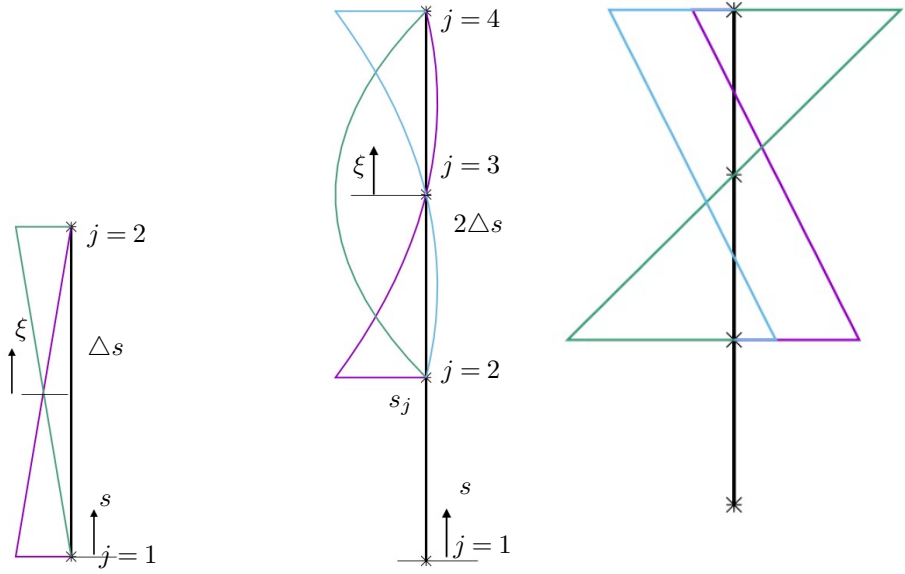

The functions $M_{j}$ are defined first as functions $\bar{M}_{j}(\xi)$ of the local coordinate $\xi$ which ranges from $-1$ to $+1$ (Fig. \ref{Parainterpol}).
For a linear interpolation the shape functions and derivatives are given by 
\begin{eqnarray}
\label{ }
\bar{M}_{1}(\xi)= 0.5(1-\xi) & \quad & \frac{\partial\bar{M}_{1}(\xi)}{\partial \xi}= -0.5 \\
\nonumber
\bar{M}_{2}(\xi)= 0.5 (1+ \xi) & \quad & \frac{\partial\bar{M}_{2}(\xi)}{\partial \xi}= 0.5
\end{eqnarray}

For a quadratic interpolation the shape function and derivatives are given by
\begin{eqnarray}
\label{ }
\bar{M}_{2}(\xi)= 1 - \xi^2  & \quad & \frac{\partial\bar{M}_{2}(\xi)}{\partial \xi}= - 2 \xi\\
\nonumber
\bar{M}_{1}(\xi)= 0.5(1-\xi) - 0.5 \ \bar{M}_{2} & \quad & \frac{\partial\bar{M}_{1}(\xi)}{\partial \xi}= -0.5(1 +  \frac{\partial\bar{M}_{2}(\xi)}{\partial \xi}) \\
\nonumber
\bar{M}_{3}(\xi)= 0.5 (1+ \xi) - 0.5 \ \bar{M}_{2} & \quad &  \frac{\partial\bar{M}_{3}(\xi)}{\partial \xi}= 0.5(1  - \frac{\partial\bar{M}_{2}(\xi)}{\partial \xi})
\end{eqnarray}

The transformation between local coordinate $s$ and $\xi$ and the associated Jacobian is for linear interpolation:
\begin{eqnarray}
s  =  \frac{1}{2} (1+ \xi) & \quad & \frac{\partial s}{\partial \xi}=\frac{1}{2}
\end{eqnarray}
and for quadratic interpolation:
\begin{eqnarray}
s  =  \frac{2 \triangle s}{2} (1+ \xi) + s_{j} & \quad & \frac{\partial s}{\partial \xi}=\triangle s
\end{eqnarray}

To compute the strains at grid points we require the derivatives of the shape functions at those points only. For example we have for the strain in local direction $s$ at grid point $i$:
\begin{equation}
\label{ }
\epsilon_{s}(s_{i})= \sum_{j=1}^{J} \frac{\partial M_{j}(s_{i})}{\partial s} u_{sj}
\end{equation}

For the quadratic interpolation, if $i$ is inside the inclusions we have the derivative to $s$:
\begin{eqnarray}
\label{ }
\frac{\partial M_{i-1}}{\partial s}&= &\frac{\partial \bar{M}_{1}(\xi=0)}{\partial \xi} \frac{\partial \xi}{\partial s}= -0.5 \frac{1}{\triangle s} \\
\frac{\partial M_{i}}{\partial s}&=& \frac{\partial \bar{M}_{2}(\xi=0)}{\partial \xi} \frac{\partial \xi}{\partial s}= 0 \\
\frac{\partial M_{i+1}}{\partial s}&=& \frac{\partial \bar{M}_{3}(\xi=0)}{\partial \xi} \frac{\partial \xi}{\partial s}= 0.5 \frac{1}{\triangle s}
\end{eqnarray}
with all other $M_{j}$ terms equal to zero.

For a point at the top we have:
\begin{eqnarray}
\label{ }
\frac{\partial M_{i}}{\partial s}&= &\frac{\partial \bar{M}_{1}(\xi=1)}{\partial \xi} \frac{\partial \xi}{\partial s}= -2 \frac{1}{\triangle s} \\
\frac{\partial M_{i-1}}{\partial s}&=& \frac{\partial \bar{M}_{2}(\xi=1)}{\partial \xi} \frac{\partial \xi}{\partial s}= 0.5\frac{1}{\triangle s} \\
\frac{\partial M_{i-2}}{\partial s}&=& \frac{\partial \bar{M}_{3}(\xi=1)}{\partial \xi} \frac{\partial \xi}{\partial s}= 1.5 \frac{1}{\triangle s}
\end{eqnarray}

For a point at the bottom we have:
\begin{eqnarray}
\label{ }
\frac{\partial M_{i-2}}{\partial s}&= &\frac{\partial \bar{M}_{1}(\xi=-1)}{\partial \xi} \frac{\partial \xi}{\partial s}= 2 \frac{1}{\triangle s} \\
\frac{\partial M_{i-1}}{\partial s}&=& \frac{\partial \bar{M}_{2}(\xi=-1)}{\partial \xi} \frac{\partial \xi}{\partial s}= -1.5\frac{1}{\triangle s} \\
\frac{\partial M_{i}}{\partial s}&=& \frac{\partial \bar{M}_{3}(\xi=-1)}{\partial \xi} \frac{\partial \xi}{\partial s}= -0.5 \frac{1}{\triangle s}
\end{eqnarray}

Gathering all strain vectors at grid points in $\{\myVecGreek{\epsilon}\}$ we can write:
\begin{equation}
\label{ }
\{\myVecGreek{\epsilon}\}= [\hat{\mathbf{ B}}] \{\myVec{u}\}
\end{equation}

After substitution of $\{\myVec{u}\}$ we obtain:
\begin{equation}
\label{Strain}
\{\myVecGreek{\epsilon}\}= [\hat{\mathbf{ B}}] \left( [\hat{\mathbf{ A}}] \mathbf{ x} + \{\bar{\mathbf{ c}}\} + [\bar{\mathbf{ B}}_{0}]\{ \myVecGreek{\sigma}_{0}\} \right)
\end{equation}

The initial stresses are computed by 
\begin{equation}
\label{ }
\myVecGreek{\sigma}_{0} =\left[ \mathbf{ D} - \mathbf{ D}_{incl}\right] \myVecGreek{\epsilon} = \left[ \mathbf{ D} - \mathbf{ D}_{incl}\right][\hat{\mathbf{ B}}] \left( [\hat{\mathbf{ A}}] \mathbf{ x} + \{\bar{\mathbf{ c}}\} + [\bar{\mathbf{ B}}_{0}]\{ \myVecGreek{\sigma}_{0}\} \right)
\end{equation}
where $\left[ \mathbf{ D} - \mathbf{ D}_{incl}\right] $ is a matrix containing $\mathbf{ D} - \mathbf{ D}_{incl}$ as sub-matrices on the diagonal.
For the interpolation function $M_{j}^{\sigma}$ we chose constant or linear basis functions similar to $M_{j}$.

\subsubsection{Computation of strain for linear inclusions}

For linear inclusions it is convenient to work with the strain in local coordinates. If we assume the bolt to be fully grouted, i.e. no slip is allowed between the bolt and the domain it is embedded in and that the Poisson's ratio of the bolt has no effect, the only strain that has to be considered is the one along the bar\footnote{It should be noted that this restriction can be lifted.}:
\begin{eqnarray}
\epsilon_{z^{\prime}} & = & \frac{\partial u_{z^{\prime}}}{\partial z^{\prime}} = \sum_{j=1}^{J} \frac{\partial M_{j}}{\partial z^{\prime}} u_{z^{\prime}j}= \sum_{j=1}^{J} \frac{\partial M_{j}}{\partial s} \frac{1}{J} ( \mathbf{ v}_{z^{\prime}}\cdot \mathbf{ u}_{j})
\end{eqnarray}
where $J$ is the Jacobian and $ \mathbf{ v}_{z^{\prime}} $ is a unit vector in $z^{\prime}$ direction.

 Eq. (\ref{Bsub}) now becomes
\begin{equation}
\label{ }
\mathbf{ B}_{j}=\frac{1}{J} \left(\begin{array}{ccc}0 & 0 & 0 \\0 & 0 & 0 \\\frac{\partial M_{j}}{\partial s}  v_{z^{\prime}_x} & \frac{\partial M_{j}}{\partial s}  v_{z^{\prime}_y}  & \frac{\partial M_{j}}{\partial s}  v_{z^{\prime}_z}  \\0 & 0 & 0 \\0 & 0 & 0 \\0 & 0 & 0\end{array}\right)
\end{equation}

The local initial stress vector is given by:
\begin{equation}
\label{ }
 \{ \myVecGreek{\sigma^{\prime}}_{0} \}= (\mathbf{ D}^{\prime} - \mathbf{ D}_{incl}^{\prime} )\{\myVecGreek{\epsilon}^{\prime} \}
\end{equation}
where
\begin{equation}
\label{ }
(\mathbf{ D}^{\prime} - \mathbf{ D}_{incl}^{\prime} )= \left(\begin{array}{cccccc}0 & 0 & 0 & 0 & 0 & 0 \\0 & 0 & 0 & 0 & 0 & 0 \\0 & 0 & E - E_{incl}  & 0 & 0 & 0 \\0 & 0 & 0 & 0 & 0 & 0 \\0 & 0 & 0 & 0 & 0 & 0 \\0 & 0 & 0 & 0 & 0 & 0\end{array}\right)
\end{equation}
where $E$ and $E_{incl}$ is the Young' modulus of the domain and the inclusion respectively and
\begin{equation}
\label{ }
 \{ \myVecGreek{\sigma^{\prime}}_{0} \}= \left\{\begin{array}{c}0 \\0 \\ \sigma_{0z^{\prime}} \\0 \\0 \\ 0 \end{array}\right\}
\end{equation}

\section{Integration of volume terms}

For the integration we have to consider 2 cases: one where point $\sourcept_{n}$ is outside the inclusion (regular integration) and one where it is not (singular integration).

\subsection{General inclusions}
The integrals to be solved are :
\begin{equation}
\mathbf{ B}_{0nj}  = \int_{\domain}  \fund{E} (\sourcept_{n},\fieldpt)  M_{j}^{\sigma}(\fieldpt) d \domain (\fieldpt) 
\end{equation}
We subdivide the inclusion into integration regions as shown in Fig. \ref{fig:Incl3D} and apply Gauss quadrature.

\subsubsection{Regular integration}
For integration region $n_{s}$ the transformation from the inclusion $\pt{s}$ coordinates to the coordinates used for Gauss integration $\bar{\myVecGreek{\xi}}=(\bar{\xi},\bar{\eta},\bar{\zeta})^{\mathrm{T}}=[-1,1]^3$ is  given by
\begin{eqnarray}
\nonumber
s & = \frac{\Delta s_{n}}{2} (1+\bar{\xi}) + s_{n_s} \\
t & = \frac{\Delta t_{n}}{2} (1+\bar{\eta}) + t_{n_s} \\
\nonumber
r& = \frac{\Delta r_{n}}{2} (1+\bar{\zeta}) + r_{n_s}
\end{eqnarray}
where $\Delta s_{n}\times \Delta t_{n} \times \Delta r_{n}$ denotes the size of the integration region and $s_{n},t_{n},r_{n}$ are the edge coordinates.
The Jacobian of this transformation is $J_{\xi}^{n}=\frac{1}{8}\ \Delta s_{n} \  \Delta t_{n} \ \Delta r_{n}$.

We can write:
\begin{equation}
  \label{GaussSecond}
   \mathbf{ B}_{0nj} = \sum_{n_{s}=1}^{N_{s}}\int_{-1}^{1} \int_{-1}^{1} \int_{-1}^{1} \fund{E} \left( \sourcept_{n},\bar{\pt{x}}(\bar{\xi},\bar{\eta}, \bar{\zeta}) \right)
 M_{j}^{\sigma} \left( \bar{\pt{x}} (\bar{\xi},\bar{\eta}, \bar{\zeta}) \right) J(\pt{s}) \ J_{\xi}^{n_{s}} \ d \bar{\xi} d \bar{\eta} d \bar{\zeta}
\end{equation}
where $ J(\mathbf{ s})$ is the Jacobian of the mapping between $\pt{s}$ and $\pt{x}$ coordinate systems.

Applying Gauss integration we have:
\begin{equation}
  \label{GaussintV}
   \mathbf{B}_{0nj}  \approx \sum_{n_{s}=1}^{N_{s}} \sum_{g_{s}=1}^{G_{s}} \sum_{g_{t}=1}^{G_{t}}  \sum_{g_{r}=1}^{G_{r}}\fund{E}\left( \sourcept_{n},\bar{\pt{x}}(\bar{\xi}_{g_{s}},\bar{\eta}_{g_{t}}, \bar{\zeta}_{g_{r}}) \right) M_{j}^{\sigma} \left( \bar{\pt{x}}(\bar{\xi}_{g_{s}},\bar{\eta}_{g_{t}}, \bar{\zeta}_{g_{r}}) \right) J(\pt{s})  \ J_{\xi}^{n_{s}} \ W_{g_{s}} \ W_{g_{t}} \ W_{g_{r}}
\end{equation}
where $N_{s}$ is the number of integration regions and $G_{s},G_{t}$ and $G_{r}$ are the number of integration points and $\bar{\xi}_{g_{s}},\bar{\eta}_{g_{t}}, \bar{\zeta}_{g_{r}}$ the Gauss point coordinates in $s, t$ and $r$ directions, respectively. $W_{g_{s}} \ W_{g_{t}} \ W_{g_{r}}$ are Gauss weights. To determine the number of Gauss points necessary for an accurate integration we consider that, whereas there is usually a moderate variation of body force, the Kernel $\fund{E}$ is $O({r}^{-2})$ so the number of integration points has to be increased if $\pt{x}_n$ is close to $\domain_0$.

\subsubsection{Singular integration}

If the integration region includes the point $\pt{x}_{n}$, then the integrand tends to infinity as the point is approached.
To deal with the integration involving the weakly singular Kernel we perform the integration in a local coordinate system, where the Jacobian tends to zero as the singularity point is approached. For this we divide the integration region into  tetrahedral sub-regions.
\begin{figure}
\begin{center}
\begin{overpic}[scale=0.5]{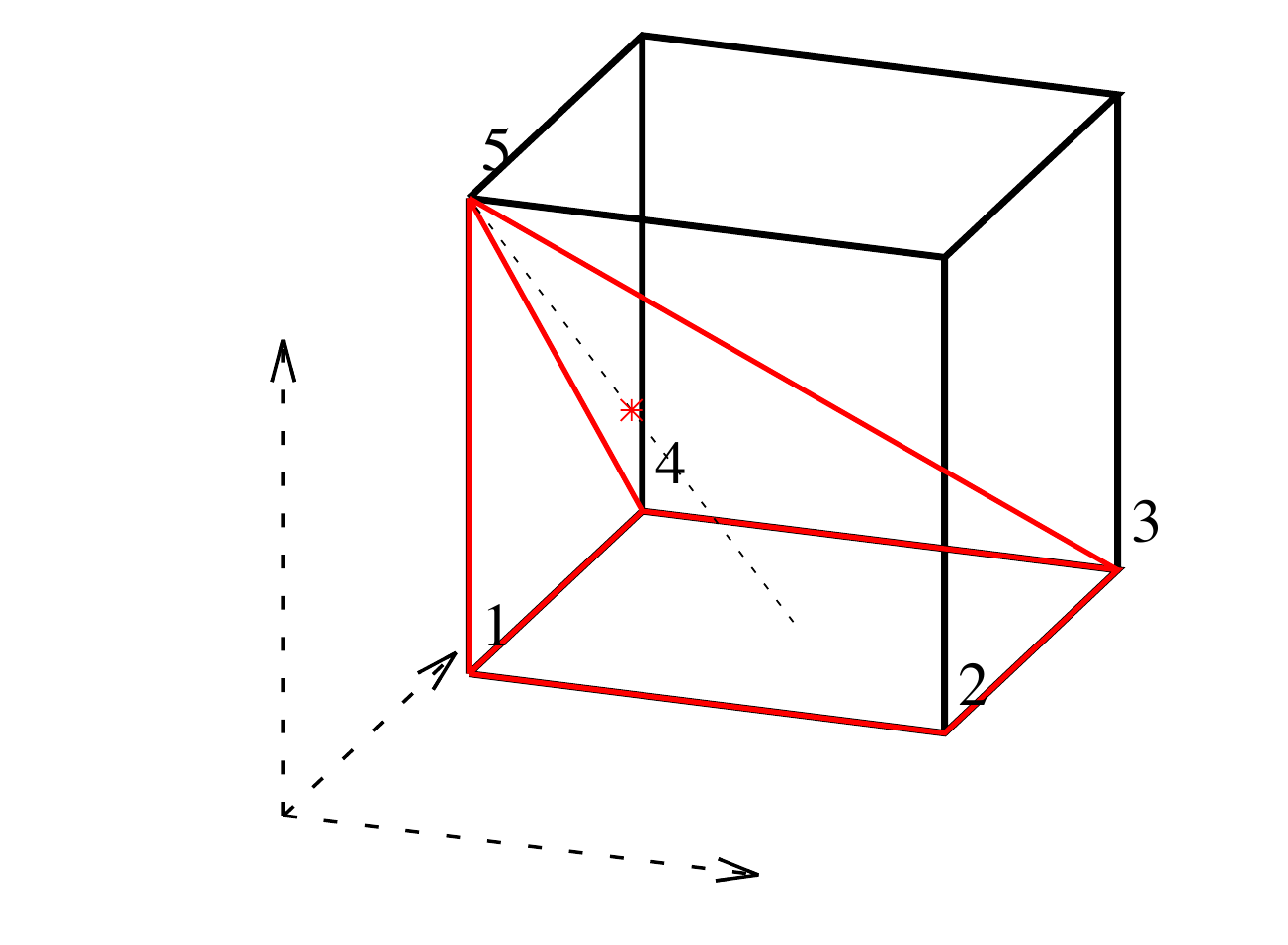}
 \put(60,5){$s$}
  \put(24,45){$r$}
    \put(30,25){$t$}
\end{overpic}
\begin{overpic}[scale=0.5]{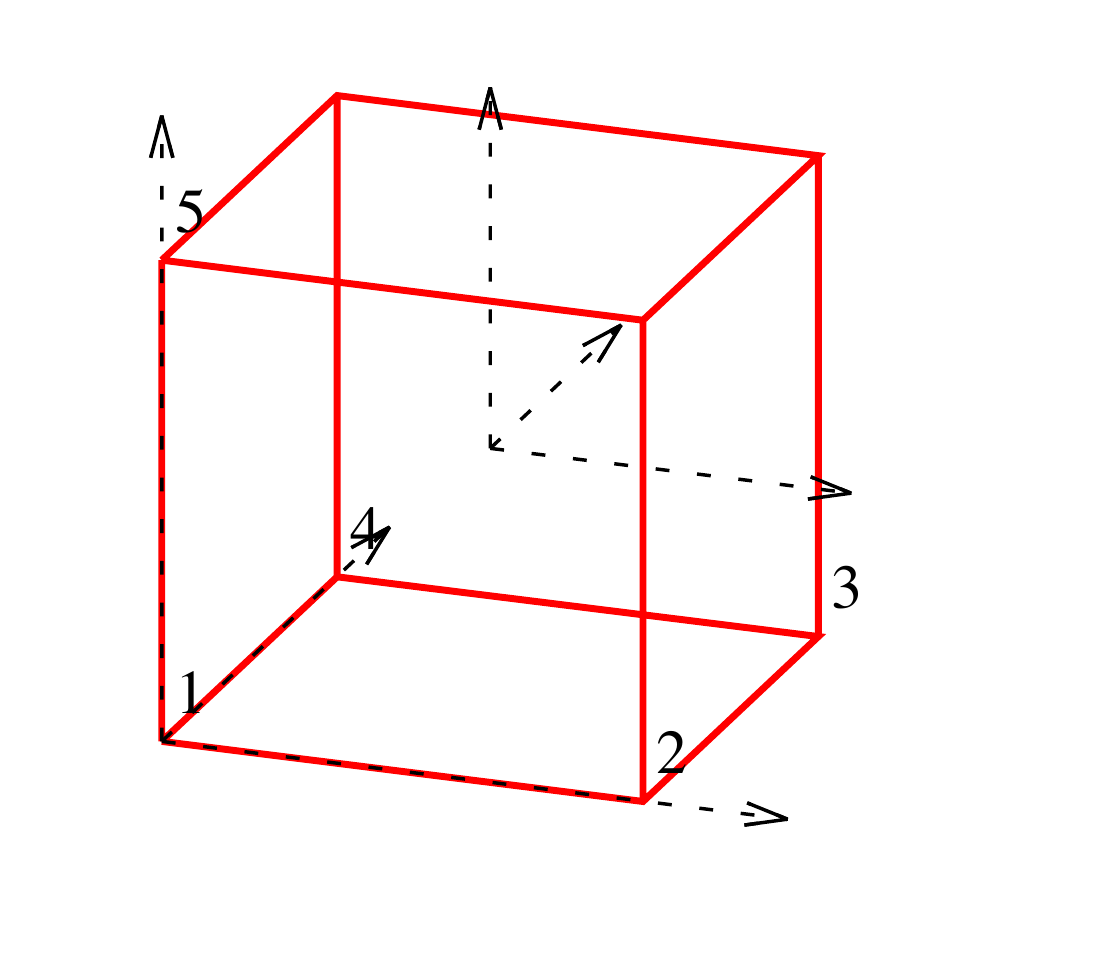}
\put(75,10){$\sigma$}
\put(10,75){$\rho$}
\put(37,37){$\tau$}
\put(80,40){$\bar{\xi}$}
\put(55,60){$\bar{\eta}$}
\put(45,75){$\bar{\zeta}$}
\end{overpic}
\caption{Singular volume integration, showing a tetrahedral subregion of an integration region and the mapping from the $\pt{s}$ to the 
$ \sigma,\tau,\rho$  coordinate system. A point with the local coordinates $\sigma=\tau=\rho=0.5$ (i.e.~$\bar{\xi}=\bar{\eta}=\bar{\zeta}=0$) is shown as a red star.}
\label{Tetra}
\end{center}
\end{figure}
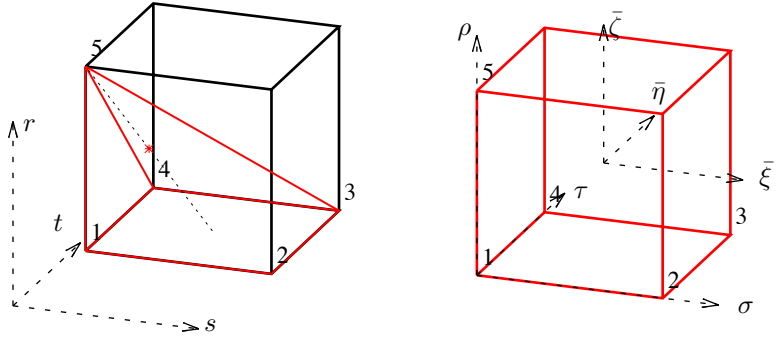
The transformation from the local $\myVecGreek{\bar{\xi}}$ coordinate system, in which the Gauss coordinates are defined, to global coordinates involves the following transformation steps:
\begin{enumerate}
  \item from $\myVecGreek{\bar{\xi}}$ to a local system $(\sigma,\tau,\rho)^{\mathrm{T}}=[0,1]^3$
  \item from ($\sigma,\tau,\rho$) to $\pt{s}$
  \item from $\pt{s}$ to $\pt{x}$
\end{enumerate}

Referring to \myfigref{Tetra} we assume that the singular point is an edge point of the integration region.
For this case the transformation is as follows:
First we determine the local coordinates $\pt{ s}_{1}$ to $\pt{ s}_{5}$ of the edge points of the tetrahedron, with 5 being the singularity point.
Next we define a linear plane NURBS surface with points 1 to 4 and map the coordinates of the point ($\sigma, \tau$) onto this surface:
\begin{equation}
%\label{ }
\pt{ s}_{0}(\sigma,\tau)= \sum_{i=1}^{4}R_{i}(\sigma,\tau) \ \pt{ s}_{i}
\end{equation}
where $R_{i}(\sigma,\tau)$ are linear basis functions.
The final map is obtained by interpolation in the $\rho$-direction:
\begin{equation}
%\label{ }
\pt{ s}(\sigma,\tau,\rho)= (1-\rho) \ \mathbf{ s}_{0}(\sigma,\tau) + \rho \  \pt{ s}_{5}
\end{equation}
The Jacobi matrix of this transformation is given by:
\begin{equation}
%\label{ }
\mathbf{ J}= \left(\begin{array}{c}(1-\rho)\frac{\partial \pt{ s}_{0}}{\partial\sigma} \\ \\ (1-\rho)\frac{\partial \pt{ s}_{0}}{\partial\tau} \\ \\ \pt{ s}_{5} - \pt{ s}_{0}\end{array}\right)
\end{equation}
The Jacobian of this transformation tends to zero as the singular point ($\rho=1$) is approached.

\subsection{Linear inclusion, reinforcement bar}
For the linear inclusions we can apply analytical integration.
In the simplest case we can assume that the initial stress is piecewise constant along the bar.
Depending on the number of internal points we divide the bar into cylindrical subregions with radius $R$ and equal length $H$ and assume the initial stress to be constant within the subdivision.
We consider two types of integration: one where point $\sourcept$ is outside the inclusion (regular integration) and one where it is not (singular integration).
Since we assume the initial stress to be constant within an integration region this means that $M_{j}^{\sigma}=1$ and the integral to be evaluated is given by:
\begin{equation}
\mathbf{ B}_{0nj}  =\int_{\domain_{j}}  \fund{E} (\sourcept_{n},\fieldpt)  d \domain_{j} (\fieldpt) 
\end{equation}
where $\domain_{j}$ denotes the subregion that corresponds with grid point $j$.

\subsubsection{Analytical computation of regular integral }

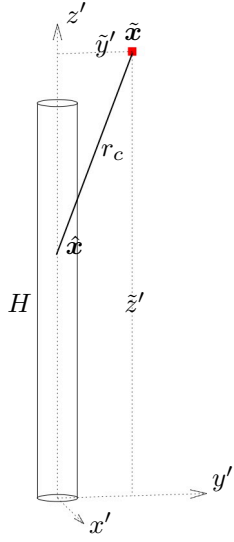
\begin{figure}
\begin{center}
\begin{overpic}[scale=0.5]{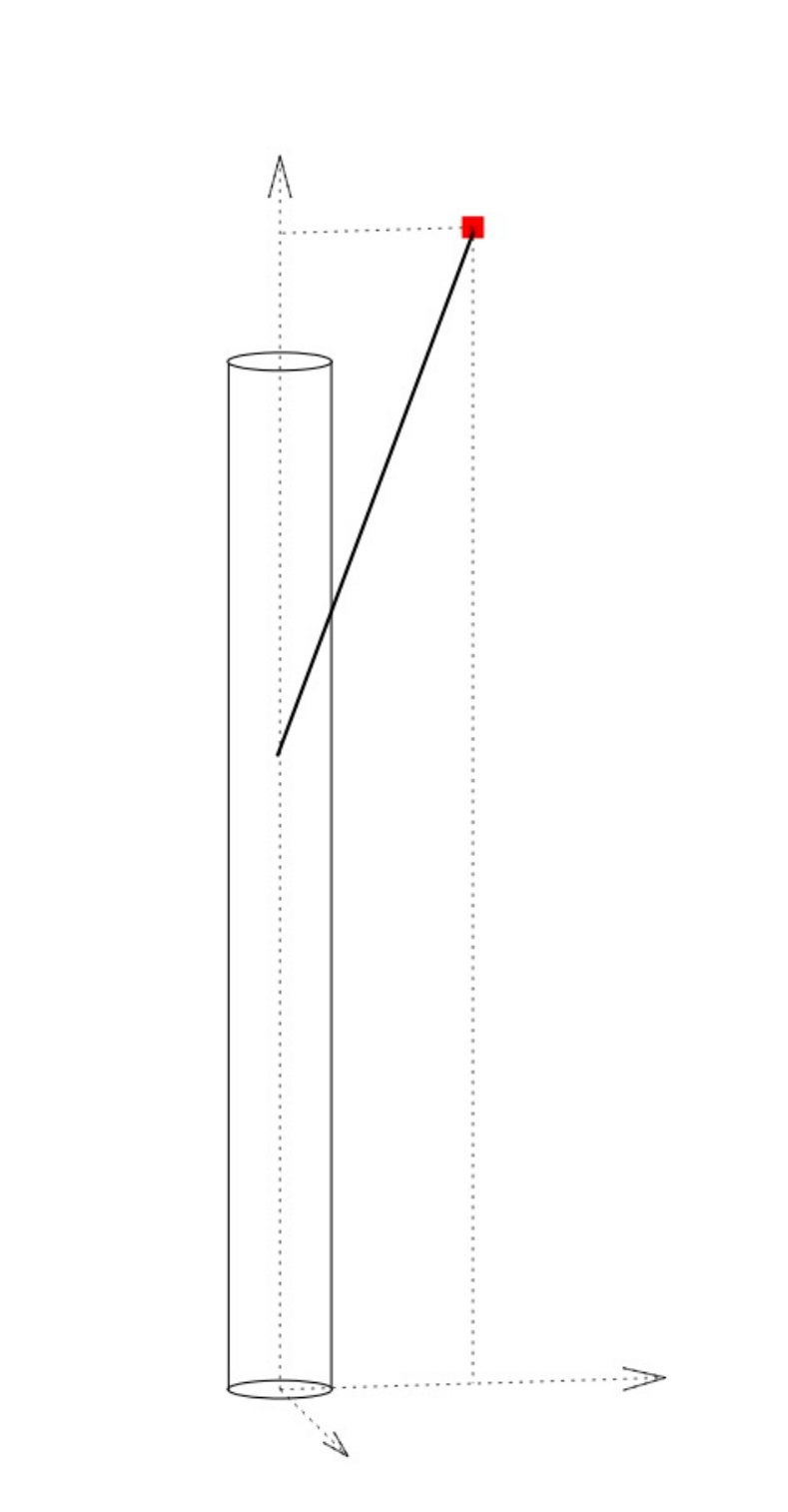}
\put(24,2){$x^{\prime}$}
\put(45,10){$y^{\prime}$}
\put(20,90){$z^{\prime}$}
\put(10,40){$H$}
\put(30,87){$\sourcept$}
\put(20,50){$\fieldpt$}
\put(30,40){$\tilde{z}^{\prime}$}
\put(25,85){$\tilde{y}^{\prime}$}
\put(26,67){$r_{c}$}
\end{overpic}
\caption{Analytical computation of regular integral for a subregion of length $H$ of a linear inclusion.}
\label{AnalR}
\end{center}
\end{figure}

Since we assume that the cross-sectional area is significantly smaller than the surrounding medium we can assume that $\fund{E}$ is constant over the cross-section. In addition we note that the result will multiply with the initial stresses in local directions ( $\{ \myVecGreek{\sigma^{\prime}}_{0} \}$). The integral to be solved is therefore:
\begin{equation}
\mathbf{ B}_{0nj}^{\prime}  = \int_{\domain_{j}} \frac{1}{r_{c}^{2}} \fund{\tilde{E}}^{\prime} d \domain_{j} (\fieldpt) 
\end{equation}
where the prime indicates that the result is computed in the local  $x^{\prime}, y^{\prime},z^{\prime}$ coordinate system (Fig. \ref{AnalR}). $r_{c}$ is the distance between the source point and a point on the axis of the inclusion and
\begin{equation}
\label{ }
\tilde{E}_{ijk}^{\prime}= -C\left[C_{3}(r_{,k}\delta_{ij} + r_{,j}\delta_{ik}) - r_{,i}\delta_{jk} + C_{4} \ r_{,i} r_{,j} r_{,k}\right]
\end{equation}

The integral to be solved is:
\begin{equation}
\label{ }
\triangle \fund{E}^{\prime}_{ij}=   -C \int_V\ \frac{1}{r_{c}^{2}} \left[C_{3}(r_{,k}\delta_{ij} + r_{,j}\delta_{ik}) - r_{,i}\delta_{jk} + C_{4} \ r_{,i} r_{,j} r_{,k}\right]   d V
\end{equation}
We choose the local axes such that $\tilde{x}^{\prime}=0$ as follows. The vector pointing in the $x^{\prime}$ direction is given by:
\begin{equation}
\label{ }
\mathbf{V}_{x^{\prime}}= (\sourcept - \fieldpt) \times \mathbf{ v}_{z^{\prime}}
\end{equation}
and the one in $y^{\prime}$ direction is:
\begin{equation}
\label{ }
\mathbf{V}_{y^{\prime}}= \mathbf{ v}_{z^{\prime}} \times \mathbf{ v}_{x^{\prime}}
\end{equation}
where the capital letter indicates that the vector is not normalised.

If point $\sourcept$ is along the axis of the bar this computation does not work and then we assume
\begin{equation}
\label{ }
\mathbf{v}_{x^{\prime}}= \mathbf{ v}_{y} \times \mathbf{ v}_{z^{\prime}}
\end{equation}
where $\mathbf{ v}_{y}$ is a vector in global $y$-direction.

We have:
\begin{equation}
\label{ }
r_1 = 0 \quad r_2 = - \tilde{y}^{\prime} \quad r_3 = z^{\prime}-\tilde{z}^{\prime} \quad r = r_c = \sqrt{\tilde{y}^{\prime 2} + (z^{\prime}-\tilde{z}^{\prime})^2}
\end{equation}
and therefore:
\begin{equation}
\label{ }
r_{,1}= 0 \quad r_{,2}= - \frac{\tilde{y}^{\prime}}{r_c} \quad r_{,3}= \frac{z^{\prime}-\tilde{z}^{\prime}}{r_c} \quad \quad dV = \pi R^2 d\tilde{z}^{\prime}
\end{equation}
The integral to be solved is:
\begin{equation}
\label{ }
\triangle \fund{E}^{\prime}_{ij}= \pi R^2 C \int_{z^{\prime}=0}^{H} \ \frac{1}{r_c^2}  \left[C_{3}(r_{,k}\delta_{ij} + r_{,j}\delta_{ik}) - r_{,i}\delta_{jk} + C_{4} \ r_{,i} r_{,j} r_{,k}\right]      dz^{\prime} 
\end{equation}

The analytical solution in Voigt notation is:
 \begin{eqnarray}
 \label{eq_DeltaEprimeReg_14}
 \triangle \fund{E^{\prime}}(1,4)(\tilde{y}^{\prime}\neq 0) & = & 2 C\pi R^2 \frac{C_3}{\tilde{y}^{\prime}}\left[  \frac{\triangle z^{\prime}}{r_{c1}} + \frac{\tilde{z}^{\prime}}{r_{c0}}  \right]   \\ 
  \nonumber
 \triangle \fund{E^{\prime}}(1,4)(\tilde{y}^{\prime}=0) & = & 0 \\
 \nonumber
 \label{eq_DeltaEprimeReg_16}
\triangle \fund{E^{\prime}}(1,6) & = & 2 C\pi R^2 C_3 \left[ \frac{1}{r_{c1}}-\frac{1}{r_{c0}}  \right]  \\
 \nonumber
 \label{eq_DeltaEprimeReg_21}
\triangle \fund{E^{\prime}}(2,1) (\tilde{y}^{\prime}\neq 0) & = & - C\pi \frac{R^2}{\tilde{y}^{\prime}} \, \left[ \frac{\triangle z^{\prime}}{r_{c1}} + \frac{z^{\prime}}{r_{c0}}  \right]  \\
 \nonumber
\triangle \fund{E^{\prime}}(2,1)(\tilde{y}^{\prime}=0) & = & 0 \\
 \nonumber
 \label{eq_DeltaEprimeReg_22}
\triangle \fund{E^{\prime}}(2,2)(\tilde{y}^{\prime}\neq 0) & = & C\pi \frac{R^2}{\tilde{y}^{\prime}} \, \left[ \left( 2(1+C_3)\tilde{y}^{\prime 2} + (1+2C_3) \triangle z^{\prime 2} \right) \frac{\triangle z^{\prime}}{r^3_{c1}} + \right. \nonumber \\ 
 \nonumber
& & \left. \frac{\tilde{z}^{\prime}}{r^3_{c0}} \left( 2(1+C_3) \tilde{y}^{\prime 2} + (1+2C_3) \tilde{z}^{\prime 2}  \right)  \right] \\
 \nonumber
\triangle \fund{E^{\prime}}(2,2)(\tilde{y}^{\prime}=0) & = & 0 
\end{eqnarray}

\begin{eqnarray}
 \label{eq_DeltaEprimeReg_23}
\triangle \fund{E^{\prime}}(2,3) & = & - C\pi R^2 \tilde{y}^{\prime} \, \left[ \frac{\triangle z^{\prime}}{r_{c1}^3} + \frac{z^{\prime}}{r_{c0}^3}  \right]  \\
 \nonumber
 \label{eq_DeltaEprimeReg_25}
\triangle \fund{E^{\prime}}(2,5) & = & 2 C\pi R^2  \, \left[ \frac{\tilde{y}^{\prime 2}+C_3 r_{c1}^2}{r_{c1}^3} - \frac{\tilde{y}^{\prime 2}+C_3 r_{c0}^2}{r_{c0}^3}  \right]  \\
 \nonumber
 \label{eq_DeltaEprimeReg_31}
\triangle \fund{E^{\prime}}(3,1) & = & C\pi R^2 \left[ \frac{1}{r_{c0}} - \frac{1}{r_{c1}}  \right]  \\
 \nonumber
 \label{eq_DeltaEprimeReg_32}
\triangle \fund{E^{\prime}}(3,2) & = & C\pi R^2 \left[ \frac{\tilde{z}^{\prime 2}}{r_{c0}^3} - \frac{\triangle z^{\prime 2}}{r_{c1}^3}  \right]  \\
 \nonumber
 \label{eq_DeltaEprimeReg_33}
\triangle \fund{E^{\prime}}(3,3) & = & C\pi R^2 \left[ \frac{(1+2C_3) \tilde{y}^{\prime 2} + 2(1+C_3)\triangle z^{\prime 2}}{r_{c1}^3} - \frac{(1+2C_3) \tilde{y}^{\prime 2} + 2(1+C_3) \tilde{z}^{\prime 2}}{r_{c0}^3} \right] \\
 \nonumber
 \label{eq_DeltaEprimeReg_35}
\triangle \fund{E^{\prime}}(3,5) (\tilde{y}^{\prime}\neq 0) & = & 2 C\pi \frac{R^2}{\tilde{y}^{\prime}} \left[ \frac{(C_3 r_{c1}^2 + \triangle z^{\prime 2}) \triangle z^{\prime}}{r_{c1}^3} + \frac{\tilde{z}^{\prime}}{r_{c0}^3} (\tilde{z}^{\prime 2}+C_3 r_{c0}^2) \right] \\
 \nonumber
\triangle \fund{E^{\prime}}(3,5)(\tilde{y}^{\prime}=0) & = & 0
\end{eqnarray}
where:
\begin{equation}
\triangle z^{\prime} = H - \tilde{z}^{\prime} \quad r_{c1} = \sqrt{\tilde{y}^{\prime 2}+\triangle z^{\prime 2}} \quad  r_{c0} = \sqrt{\tilde{y}^{\prime 2}+z^{\prime 2}}
\end{equation}

Since the result of the multiplication with $\myVecGreek{\sigma^{\prime}}_{0}$ has to be in global coordinates a transformation to the global system is necessary:
\begin{equation}
\label{ }
\mathbf{ B}_{0nj}= \mathbf{ T} \mathbf{ B}_{0nj}^{\prime}
\end{equation}
where $\mathbf{ T}$ is the transformation matrix given by:
\begin{equation}
\label{ }
\mathbf{ T}= \left(\begin{array}{ccc}v_{x^{\prime}_{x}} & v_{y^{\prime}_{x}} & v_{z^{\prime}_{x}} \\v_{x^{\prime}_{y}} & v_{y^{\prime}_{y}} & v_{z^{\prime}_{y}} \\v_{x^{\prime}_{z}} & v_{y^{\prime}_{z}} & v_{z^{\prime}_{z}}\end{array}\right)
\end{equation}

\subsubsection{Analytical computation of singular integral}
\begin{figure}
\begin{center}
\begin{overpic}[scale=0.5]{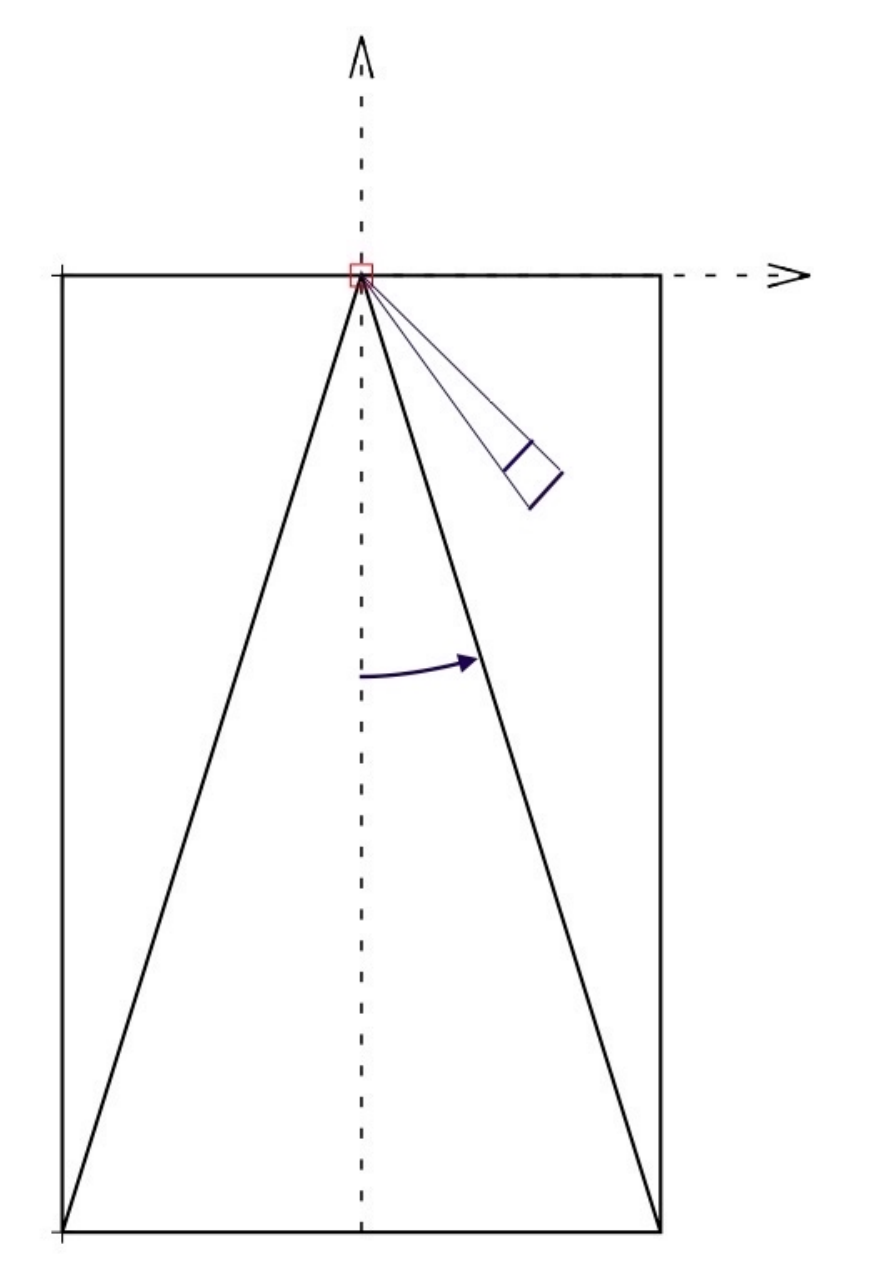}
\put(36,72){$r$}
\put(43,65){$dr$}
\put(31,55){$\theta$}
\put(40,5){$R$}
\put(44,60){$ r d \theta $}
\put(30,30){1}
\put(10,50){2}
\put(55,40){$H$}
\put(31,97){$z^{\prime}$}
\put(60,85){$y^{\prime}$}
\end{overpic}
\begin{overpic}[scale=0.4]{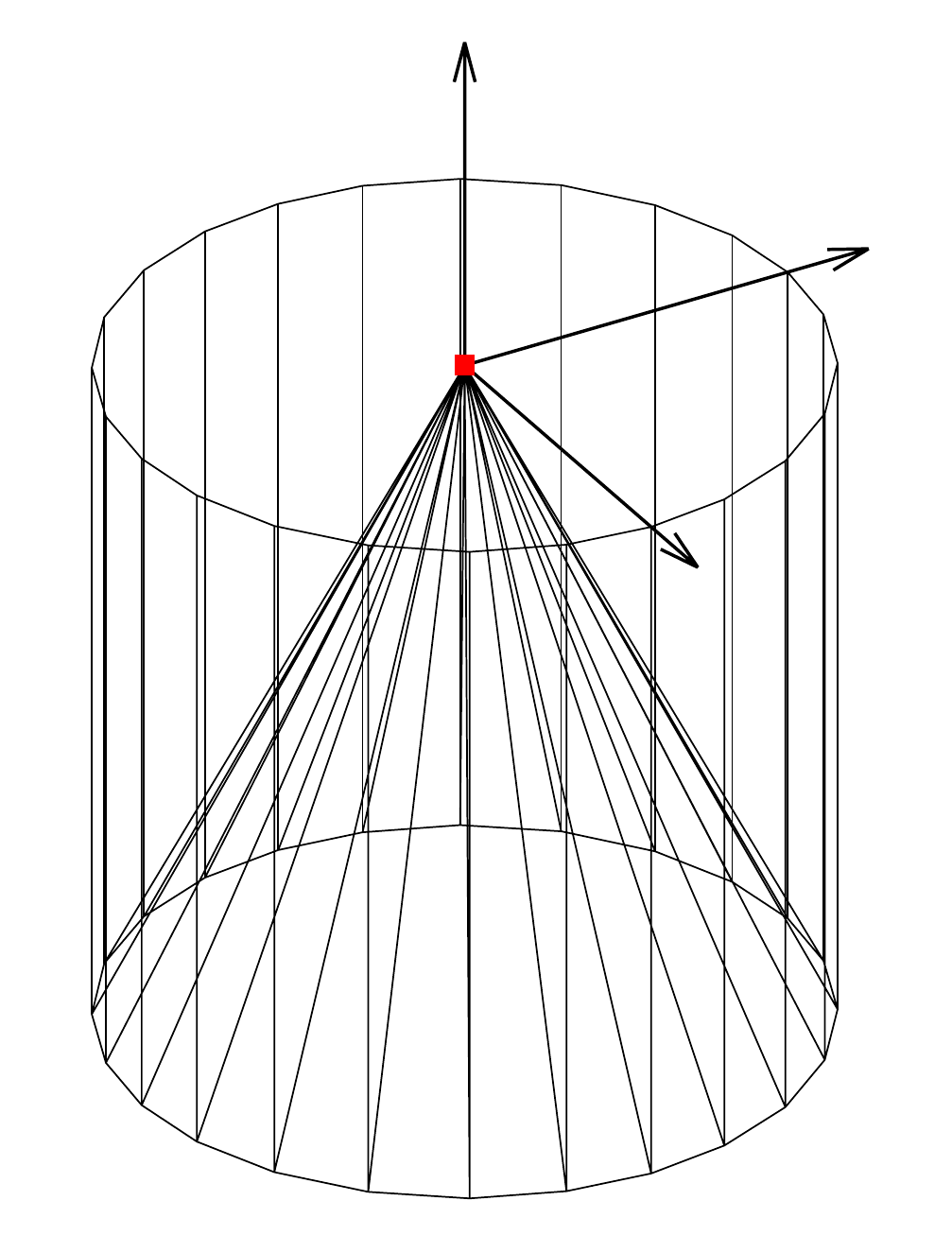}
\put(31,97){$z^{\prime}$}
\put(60,85){$y^{\prime}$}
\put(50,50){$x^{\prime}$}
\end{overpic}
\caption{Explanation of singular integration by subdivision into conical subregions. Left: section through bar, right: Axonometric view. Singular point is marked by red square}
\label{Sing2}
\end{center}
\end{figure}
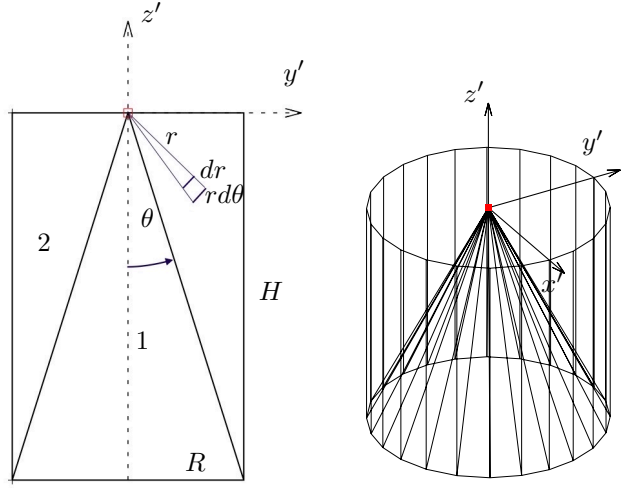

Referring to Fig. \ref{Sing2} we subdivide the bolt into 2 subregions and obtain the following integrals in polar coordinates:
\begin{equation}
\label{I1I2}
\triangle \fund{E^{\prime}}_1  =  \int_{\phi=0}^{2\pi} \int_{\theta=\tilde{\theta}}^{0} \int_{r=0}^{\frac{H}{\cos{\theta}}}  \frac{1}{r^{2}} \fund{\tilde{E}^{\prime}} \sin{\theta}  dr \ r^{2} d\theta d\phi 
\end{equation}
\begin{equation}
\triangle \fund{E^{\prime}}_2  =  \int_{\phi=0}^{2\pi} \int_{\theta=\pi/2}^{\tilde{\theta}} \int_{r=0}^{\frac{R}{\sin{\theta}}} \frac{1}{r^{2}}  \fund{\tilde{E}^{\prime}} \sin{\theta}  dr \ r^{2} d\theta d\phi \nonumber
\end{equation}
with $\tilde{\theta}= \arctan (R/H)$. It can be seen that the $r^{2}$ terms cancel out which means that the integrand is no longer singular.

 The terms  of $\triangle \fund{E}^{\prime} = \triangle \fund{E}^{\prime}_1+\triangle \fund{E}^{\prime}_2$, in Voigt notation, different from zero are:
 \begin{eqnarray}
 \label{eq_DeltaEprime_1625}
 \triangle \fund{E^{\prime}}(1,6) = \triangle \fund{E^{\prime}}(2,5) & = & \frac{C\pi}{2} \left[ H \left( 8+8 C_3-(9+8 C_3) \cos{\tilde{\theta}} +  \cos{(3\tilde{\theta})} \right) - \right. \nonumber \\
 & & \left. - 4 R \left( -1 -2 C_3 + 2 C_3 \sin{\tilde{\theta}} + \sin^3{\tilde{\theta}} \right) \right] \\
  \nonumber
 \label{eq_DeltaEprime_3132}
 \triangle \fund{E^{\prime}}(3,1) = \triangle \fund{E^{\prime}}(3,2) & = & -C\pi \left[ R + \sin{\tilde{\theta}} \left( \frac{H}{2} \sin{(2\tilde{\theta})} + R (\sin^2{\tilde{\theta}}-2) \right) \right]  \\
  \nonumber
 \label{eq_DeltaEprime_33}
 \triangle \fund{E^{\prime}}(3,3) & = & C\pi \left[ -2 H (\cos{\tilde{\theta}}-1) (2 C_3 +
\cos{\tilde{\theta}}+\cos^2{\tilde{\theta}}) \right. \nonumber \\
 \nonumber
& & \left. + R (2+4 C_3 - (3+4 C_3 +\cos{(2\tilde{\theta})}) \sin{\tilde{\theta}}) \right]
\end{eqnarray}
As before a transformation to the global system is necessary:
\begin{equation}
\label{ }
\mathbf{ B}_{0nj}= \mathbf{ T} \mathbf{ B}_{0nj}^{\prime}
\end{equation}

\section{Solution procedure}
Eqs. (\ref{DisIE}), (\ref{Strain}) form a linear system of equations where the initial stresses are a function of the elastic strains. The system of equations may be solved iteratively using a modified Newton-Raphson method. If the inclusions are elastic, however, it is possible solve the system in one step as will be shown.

\subsection{Iterative solution using modified Newton Raphson}
Using a modified Newton-Raphson method first solve
\begin{equation}
\label{DisIE1}
[\mathbf{ L}] \{\mathbf{ x}\}_{0} = \{\mathbf{ r}\} 
\end{equation}
and then compute increments of the boundary unknown $\{\mathbf{ x}\}_{i}$ using
\begin{equation}
\label{DisIE2}
[\mathbf{ L}] \{\mathbf{ x}\}_{i} = [\mathbf{ B}_{0}] \{ \myVecGreek{\sigma}_{0}\} 
\end{equation}
where the subscript $i$ is the iteration number.
The final values are obtained by summing all the increments:
\begin{equation}
\label{Sumiter}
\{\mathbf{ x}\}_{i} = \{\mathbf{ x}\}_{0} + \{\mathbf{ x}\}_{1} + \{\mathbf{ x}\}_{2} \cdots
\end{equation}

\subsection{One step solution}
A one step solution is possible by combining equation (\ref{DisIE}) with (\ref{Strain}).
Eq. (\ref{Strain}) can be written in the following form:
\begin{equation}
\label{Strain_rev1}
\{\myVecGreek{\epsilon}\}= [\hat{\mathbf{ C}}]  \{ \mathbf{ x} \} + \{\bar{\bar{\mathbf{ c}}}\} + [\bar{\mathbf{ C}}_{0}]  ( [\mathbf{ D}] - [\mathbf{ D}_{incl}]) \{ \myVecGreek{\epsilon} \}
\end{equation}
\noindent where:
\begin{equation}
\label{}
[\hat{\mathbf{ C}}] = [\hat{\mathbf{ B}}] [\hat{\mathbf{ A}}] \hspace{10mm}
[\hat{\mathbf{ C}}_{0}] = [\hat{\mathbf{ B}}] [\bar{\mathbf{B}}_{0}] \hspace{10mm}
\{\bar{\bar{\mathbf{ c}}}\} = [\hat{\mathbf{ B}}] \{\bar{\mathbf{ c}}\}
\end{equation}
Eq. (\ref{Strain_rev1}) along with Eq. (\ref{DisIE}) form the following linear system of equations:
\begin{equation}
\label{Final_sys}
\begin{pmatrix}
  [\mathbf{ L}]        & - [\mathbf{ B}_{0}] ( [\mathbf{ D}] - [\mathbf{ D}_{incl}])  \\ \\
 - [\hat{\mathbf{ C}}] & [{\mathbf{I}}] - [\hat{\mathbf{ C}}_{0}] ( [\mathbf{ D}] - [\mathbf{ D}_{incl}])
\end{pmatrix}
\begin{pmatrix}
\{\mathbf{ x}\} \\ \\
\{\myVecGreek{\epsilon}\}
\end{pmatrix}
=
\begin{pmatrix}
\{\mathbf{ r}\} \\ \\
\{\bar{\bar{\mathbf{ c}}}\}
\end{pmatrix}
\end{equation}
that can be solved in terms of boundary unknowns and internal strains.

However, it is possible to obtain a system of equations that only multiplies with the boundary unknown:
\begin{equation}
\label{onestep}
[\mathbf{ L}]^{\prime} \{\mathbf{ x}\} = \{\mathbf{ r}\}^{\prime} 
\end{equation}
where $[\mathbf{ L}]^{\prime}$ and $\{\mathbf{ r}\}^{\prime} $ are modified left and right hand sides.

We rewrite the strain vector as: 
\begin{equation}
\label{Strain_rev2}
\{\myVecGreek{\epsilon}\}= ([\mathbf{I}] - 
[\hat{\mathbf{ C}}_0] ([\mathbf{D}]-[\mathbf{D}_{incl}])^{-1}) ([\hat{\mathbf{C}}] \{\mathbf{ x}\} + \{\bar{\bar{\mathbf{ c}}}\}) = [\mathbf{A}] \{\mathbf{x}\} + \{\mathbf{b}\}
\end{equation}

\noindent where

\begin{equation}
[\mathbf{A}] = ( [\mathbf{I}] - [\hat{\mathbf{C}}_0] ([\mathbf{D}] - [\mathbf{D}_{incl}]))^{-1} [\hat{\mathbf{C}}] \hspace{5mm} \{\mathbf{b}\} = ( [\mathbf{I}] - [\hat{\mathbf{C}}_0] 
([\mathbf{D}] - [\mathbf{D}_{incl}]))^{-1} \{\bar{\bar{\mathbf{ c}}}\}
\end{equation}

Eq. (\ref{Strain_rev2}) can be inserted in Eq. (\ref{DisIE}) in order to obtain:
\begin{equation}
[\mathbf{L}] \{\mathbf{x}\} = \{\mathbf{r}\} + [\mathbf{B}_0] 
([\mathbf{D}] - [\mathbf{D}_{incl}]) ([\mathbf{A}] \{\mathbf{x}\} + \{\mathbf{b}\})
\end{equation}
and, hence, the following system of equations can be obtained
\begin{equation}
( [\mathbf{L}] - [\mathbf{B}_0]) ([\mathbf{D}] - [\mathbf{D}_{incl}]) [\mathbf{A}] ) 
\{\mathbf{x}\} = \{\mathbf{r}\} + [\mathbf{B}_0] 
([\mathbf{D}] - [\mathbf{D}_{incl}]) \{\mathbf{b}\}
\end{equation}

The matrices in Eq. (\ref{onestep}) are defined by:
\begin{eqnarray}
[\mathbf{L}]^{\prime} & = & ( [\mathbf{L}] - [\mathbf{B}_0]) ([\mathbf{D}] - [\mathbf{D}_{incl}]) [\mathbf{A}] )   \\
\{\mathbf{ r}\}^{\prime} & = &\{\mathbf{r}\} + [\mathbf{B}_0] 
([\mathbf{D}] - [\mathbf{D}_{incl}]) \{\mathbf{b}\} 
\end{eqnarray}

\begin{remark}
The modified left hand side $[\mathbf{L}]^{\prime}$ can be substituted for $[\mathbf{L}]$ in Eq. (\ref{DisIE1}) and a Newton-Raphson iteration applied to the case where the inclusions exhibit elasto-plastic behaviour in addition to having elastic material properties that are different to the domain. So the idea is to first apply a one step solution to account for the difference in elastic properties and then solve the nonlinear problem. Indeed, there is even a possibility that the left hand side is modified using the elasto-plastic constitutive matrix $\mathbf{ D}_{ep,incl}$ via Eq. (\ref{InitialSp}) during the iterations resulting in a true Newton-Raphson approach for the non-linear problem.
\end{remark}

\newpage
\section{Test Example 1}
The first test example is designed to test the influence of one bar on the deformation of a cube.
\begin{figure}[h]
\begin{center}
\begin{overpic}[scale=0.4]{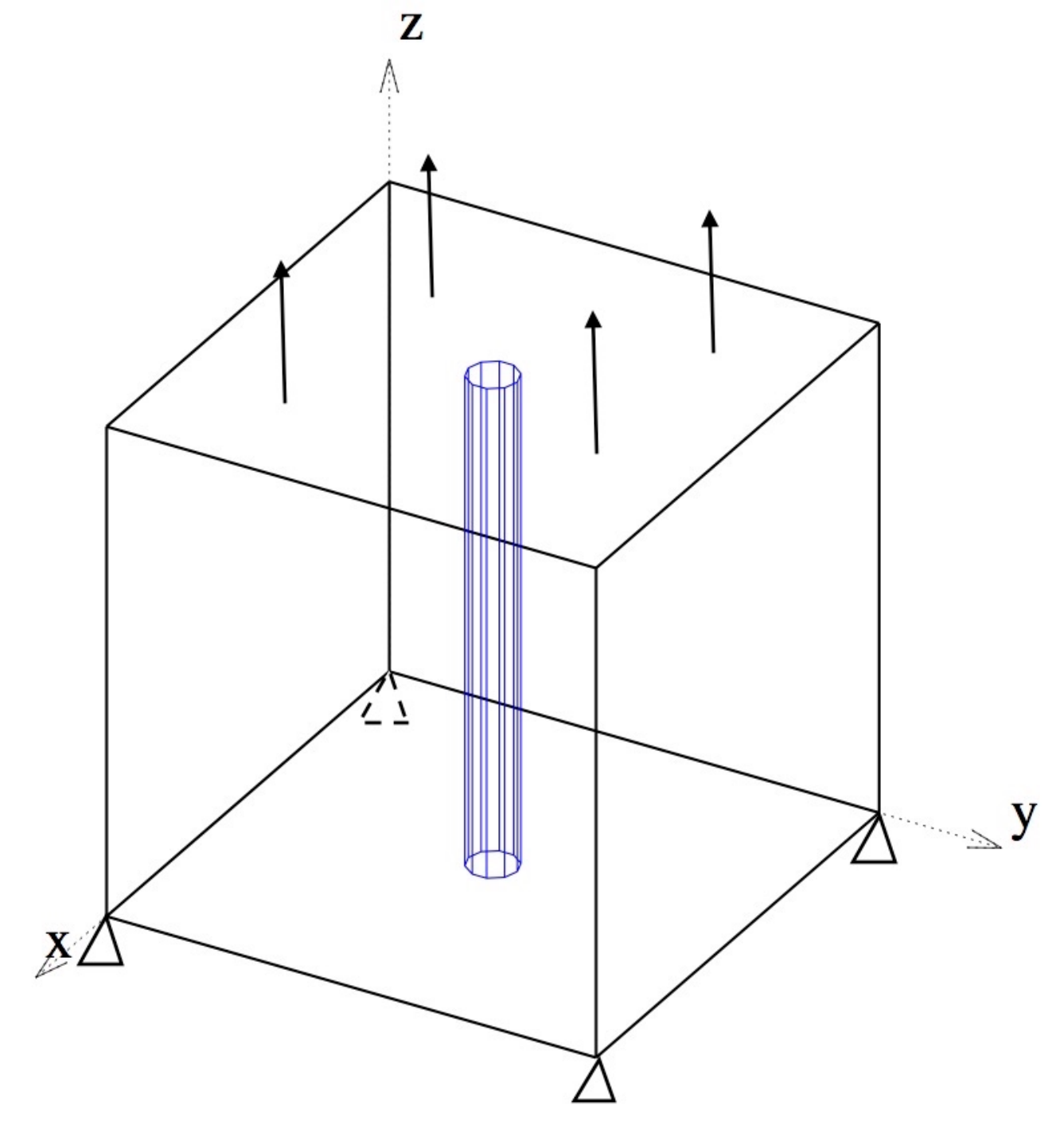}
\put(45,83){$t_{z}=1$}
\end{overpic}
\begin{overpic}[scale=0.4]{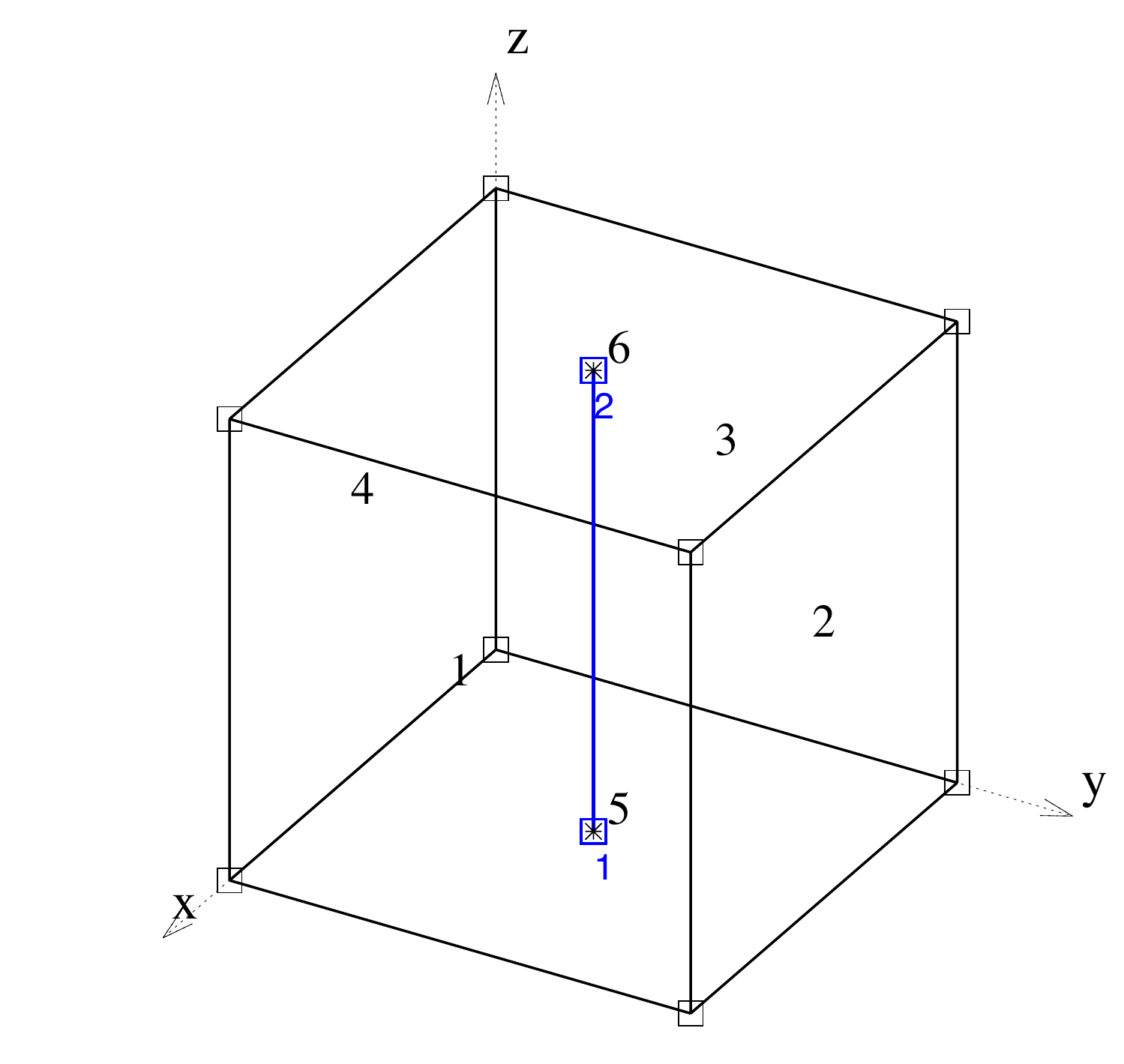}
\end{overpic}
\begin{overpic}[scale=0.4]{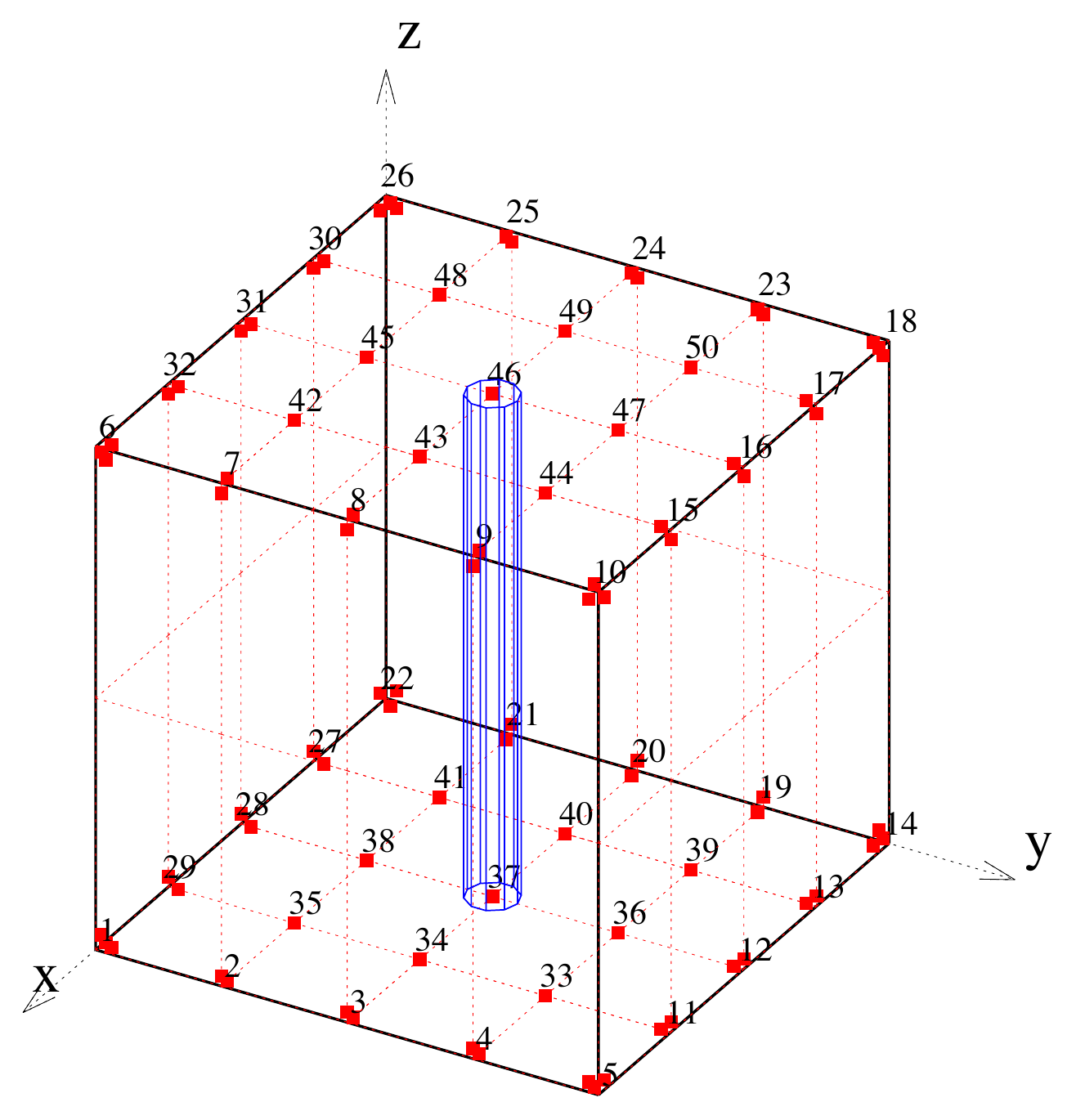}
\end{overpic}
\caption{Test example 1: Basic setup (top left) and (top right) discretisation into 6 patches, showing control points as hollow squares, and one  linear inclusion with 2 internal points. Bottom figure shows the location of collocation points after refinement.}
\label{Ex1geo}
\end{center}
\end{figure}

The basic set-up is shown in Fig. \ref{Ex1geo}.  It consists of a cube of dimension $1\times1\times1$ which is fixed at the bottom and subjected to a tensile load of 1 at the top\footnote{We use dimensionless units for this example.}.
The Youngs modulus of the cube is 1 and the Poisson's ratio 0. At the centre of the cube there is a bar with a cross-sectional diameter of 0.1 and a Young's modulus of 2.

The discretisation into 6 linear patches is shown in Fig. \ref{Ex1geo}. The inclusion was defined as a straight line with 2 control points. The number of internal points and therefore the number of integration regions along the bar was varied from 2 to 21. Constant variation of the initial stress and analytical integration was assumed within an integration region for the linear inclusion.
\begin{figure}
\begin{center}
\includegraphics[scale=0.17]{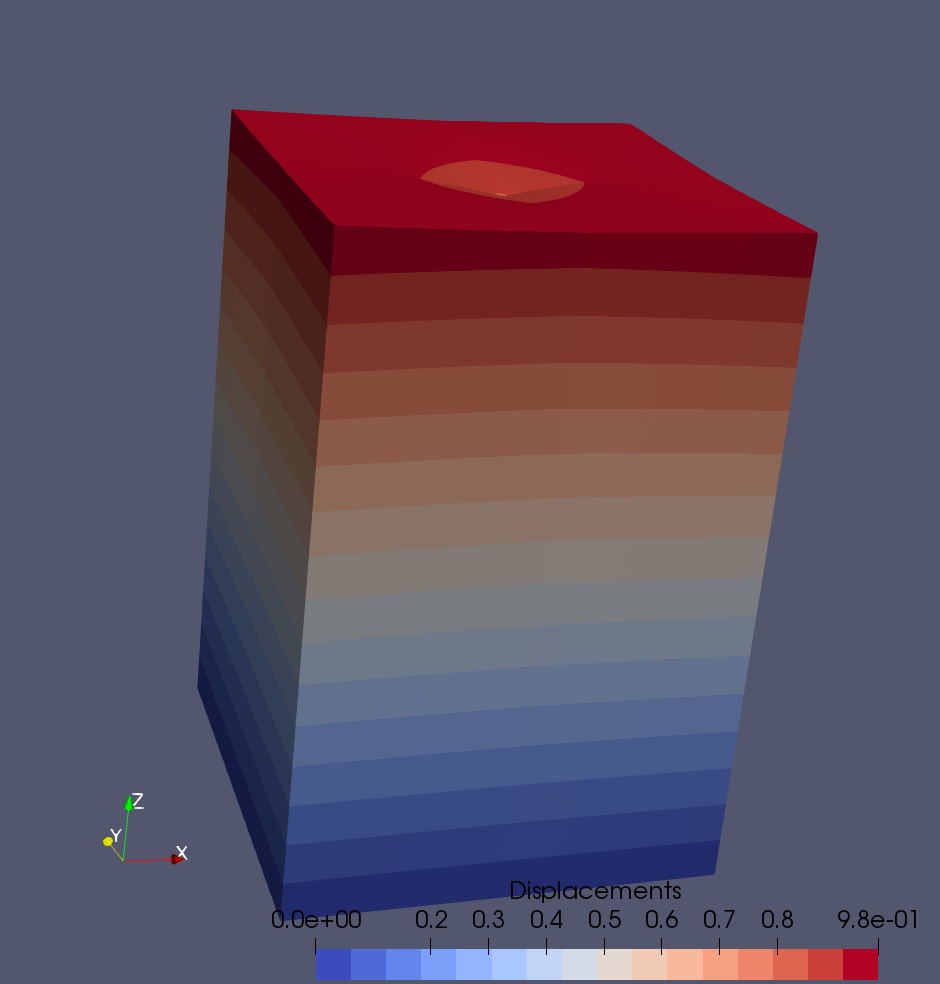}
\caption{Example 1: Displaced shape with 21 internal points for the bar.}
\label{Ex1analytical}
\end{center}
\end{figure}

For the simulation the basis functions for describing the boundary displacements were defined as follows:
The linear basis functions describing the geometry of the cube were oder elevated by one order (from linear to quadratic). At the location where the linear inclusion touches the boundary surface a knot was inserted to reduce the continuity from $C^{1}$ to $C^{0}$. The resulting location of the collocation points is shown in Fig. \ref{Ex1geo} on the bottom. The discretisation has 150 degrees of freedom.

\subsubsection{Results}
For this example we have chosen the one step solution and the resulting displaced shape is shown in Fig. \ref{Ex1analytical}.
\begin{figure}[h]
\begin{center}
\includegraphics[scale=0.6]{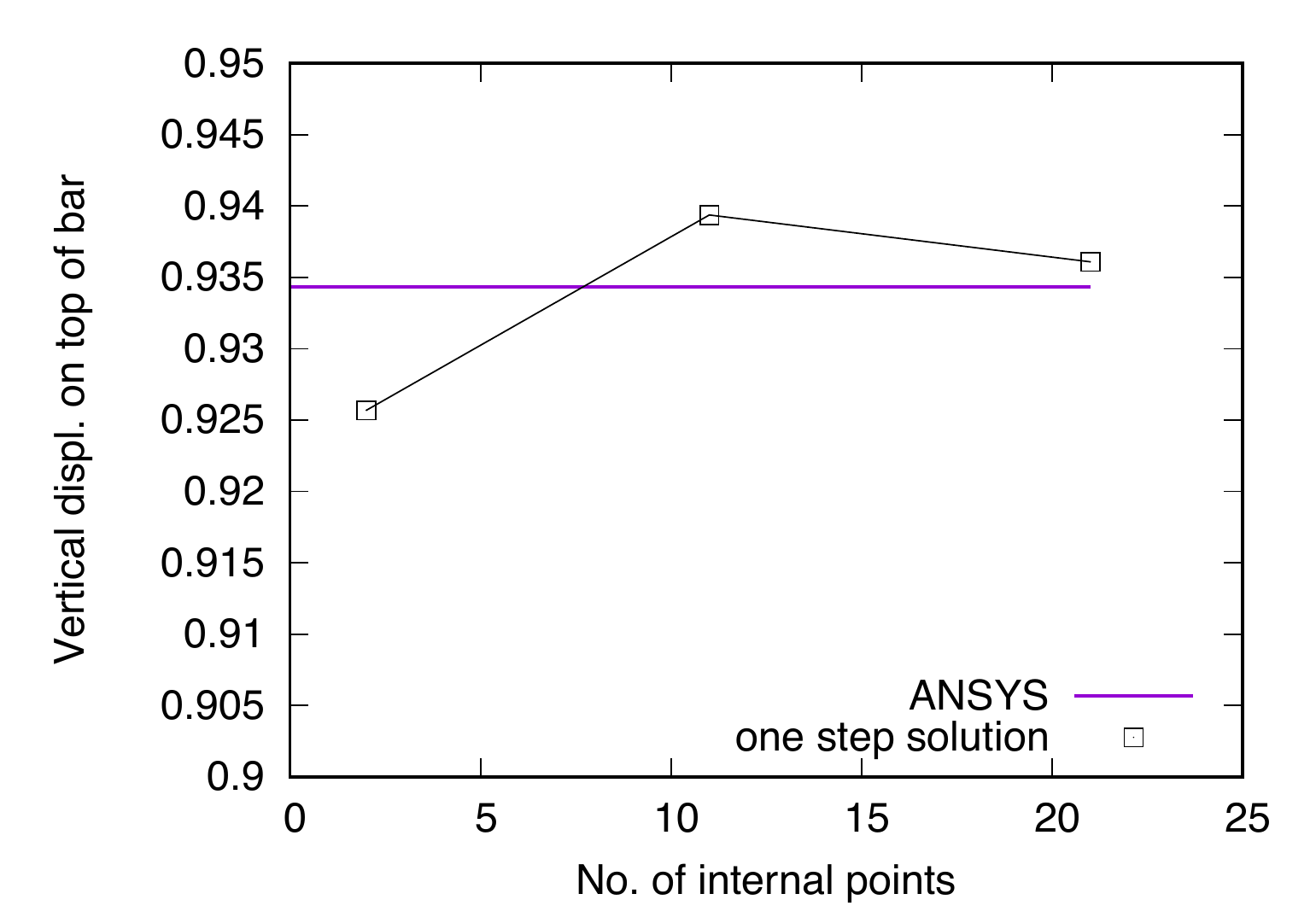}
\caption{Change of vertical displacement at the top of the bar depending on the number of internal points.}
\label{Ex1conv}
\end{center}
\end{figure}
Next we examine the effect of the number of internal points inside the linear inclusion on the results. This example, where the main variation of the displacement occurs near the top of the bar, is particularly sensitive to this parameter. Note that the number of internal points is linked to the number of integration regions over which we assume the initial stresses to be constant,
In Fig. \ref{Ex1conv} we show the convergence of the displacement towards the result of a Finite Element analysis with 19530 degrees of freedom using the ANSYS software. 

\section{Test Example 2}
\begin{figure}
\begin{center}
\includegraphics[scale=0.45]{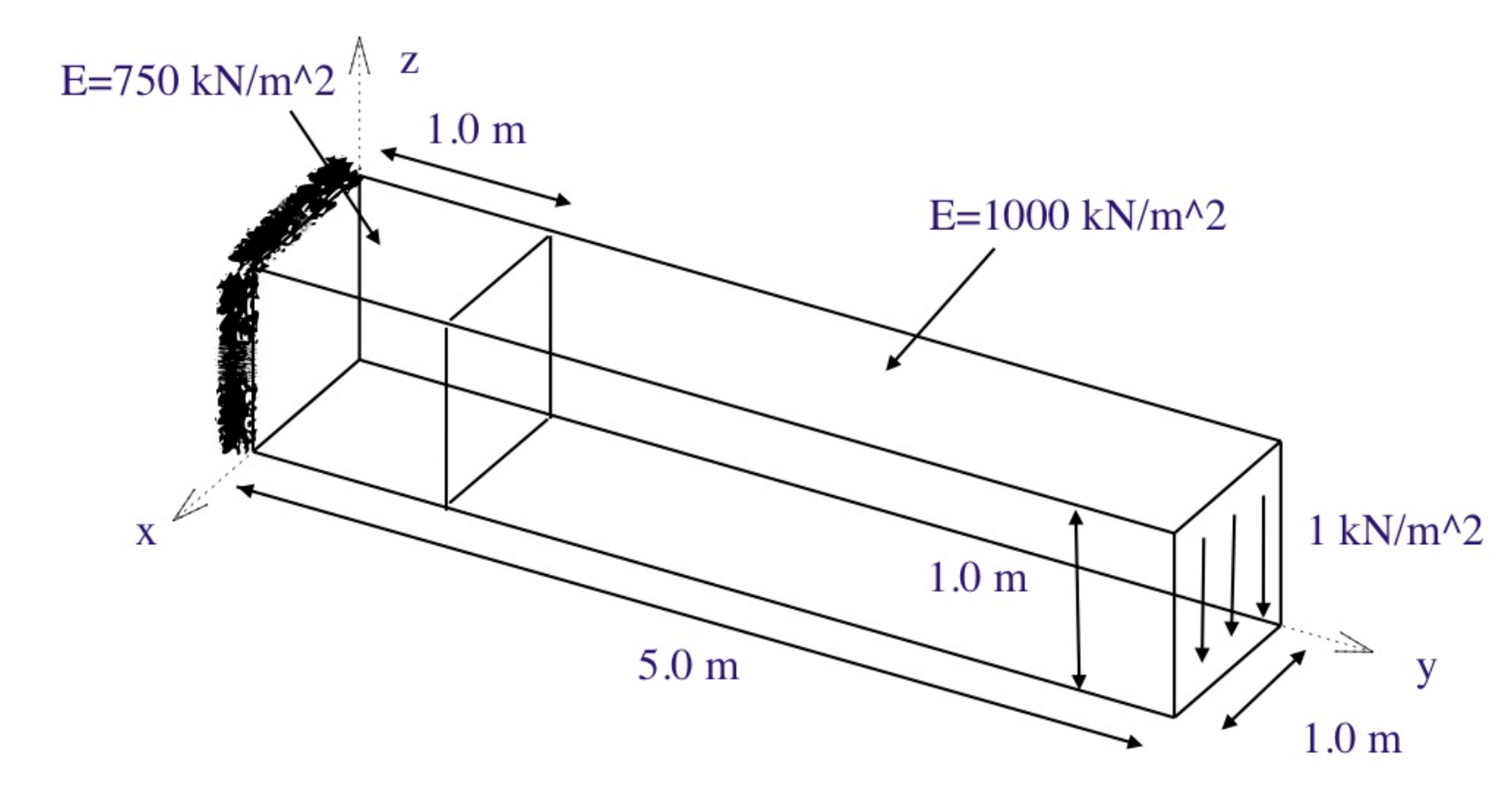}
\includegraphics[scale=0.45]{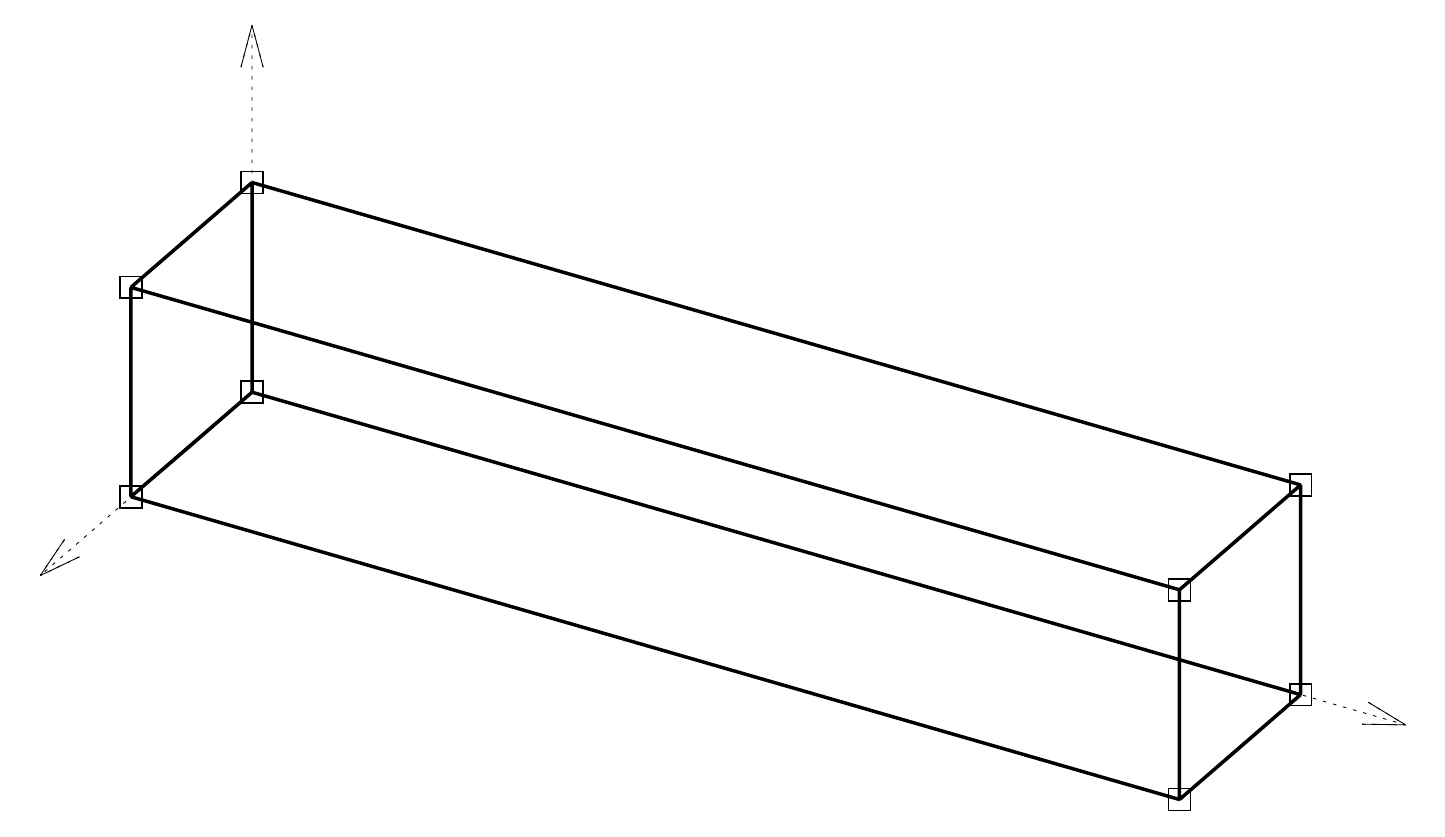}
\caption{Left:Geometry and boundary conditions of cantilever beam. Right: Geometry discretisation of the problem into 6 linear patches.}
\label{Ex2geo}
\end{center}
\end{figure}

This test example is designed to compare the modified Newton-Raphson solution with the one step solution for the case where two elastic materials exist. 
The example is a cantilever beam that consists of 2 different materials as shown in Fig. \ref{Ex2geo} and is loaded at the end by a distributed load. The Poisson's ratio is assumed to be zero. 

The geometry of the problem is defined by 6 linear NURBS patches as shown in Fig. \ref{Ex2geo}.
\begin{figure}
\begin{center}
\includegraphics[scale=0.5]{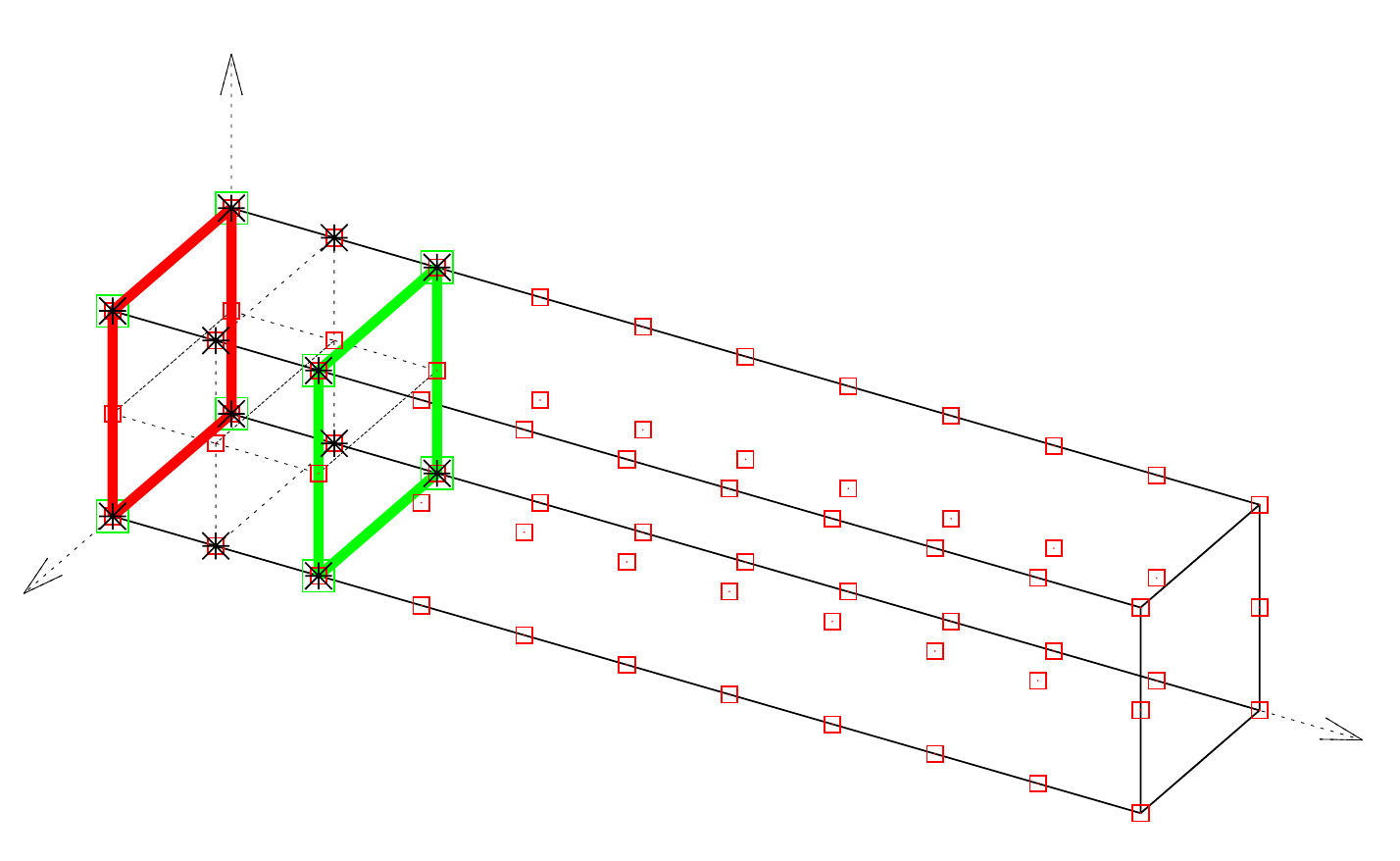}
\caption{Definition of the inclusion with reduced modulus showing "top" and "bottom" bounding surfaces and internal points. Also shown are the location of the collocation points as red squares.}
\label{Ex2patch}
\end{center}
\end{figure}
For the simulation the basis functions used for describing the variation of the boundary displacements were refined as follows:
In the directions along the cantilever the order was elevated from linear to quadratic and 4 knots were inserted.
In the vertical direction the order was elevated from linear to quadratic. The simulation has 196 degrees of freedom and the resulting collocation point locations are marked as red squares in Fig. \ref{Ex2patch}.
\begin{figure}
\begin{center}
\includegraphics[scale=0.2]{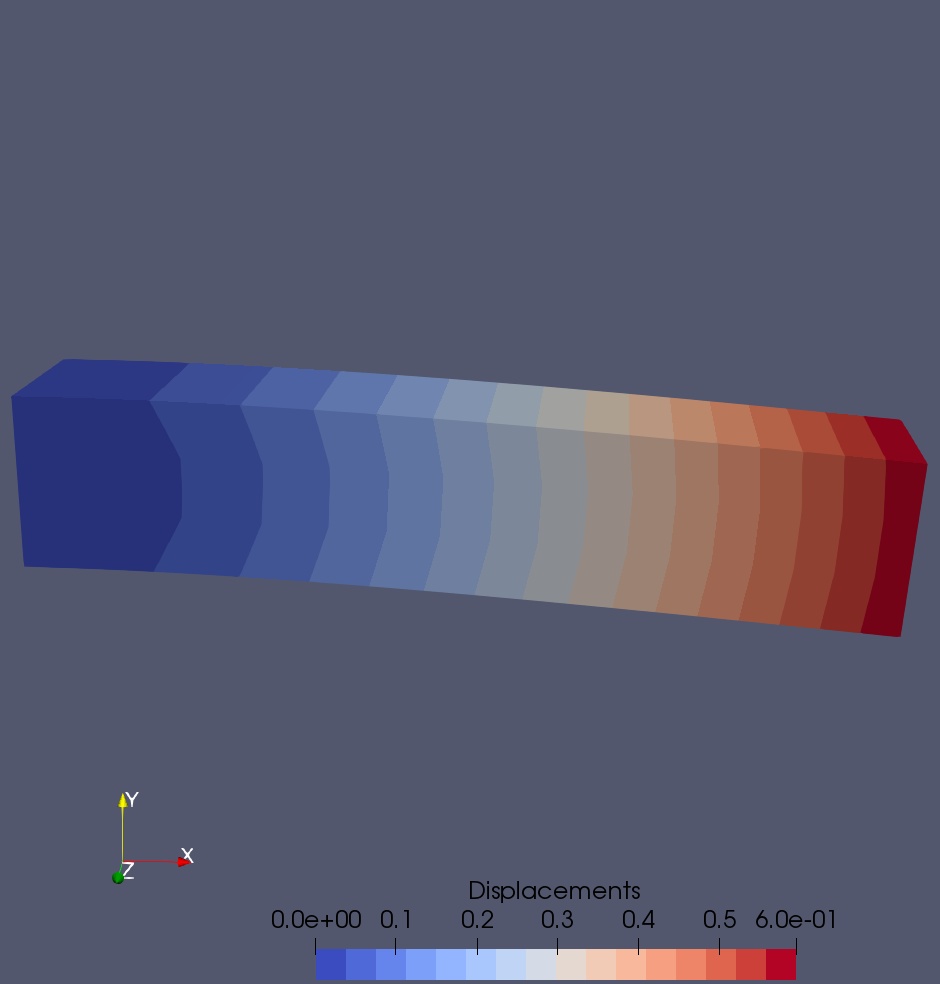}
\caption{Resulting displaced shape.}
\label{Ex2displ}
\end{center}
\end{figure}
For the fundamental solution E=1000 was used and the domain with the reduced modulus was defined as an inclusion defined by 2 NURBS surfaces as shown in Fig. \ref{Ex2patch}.

\begin{figure}
\begin{center}
\includegraphics[scale=0.5]{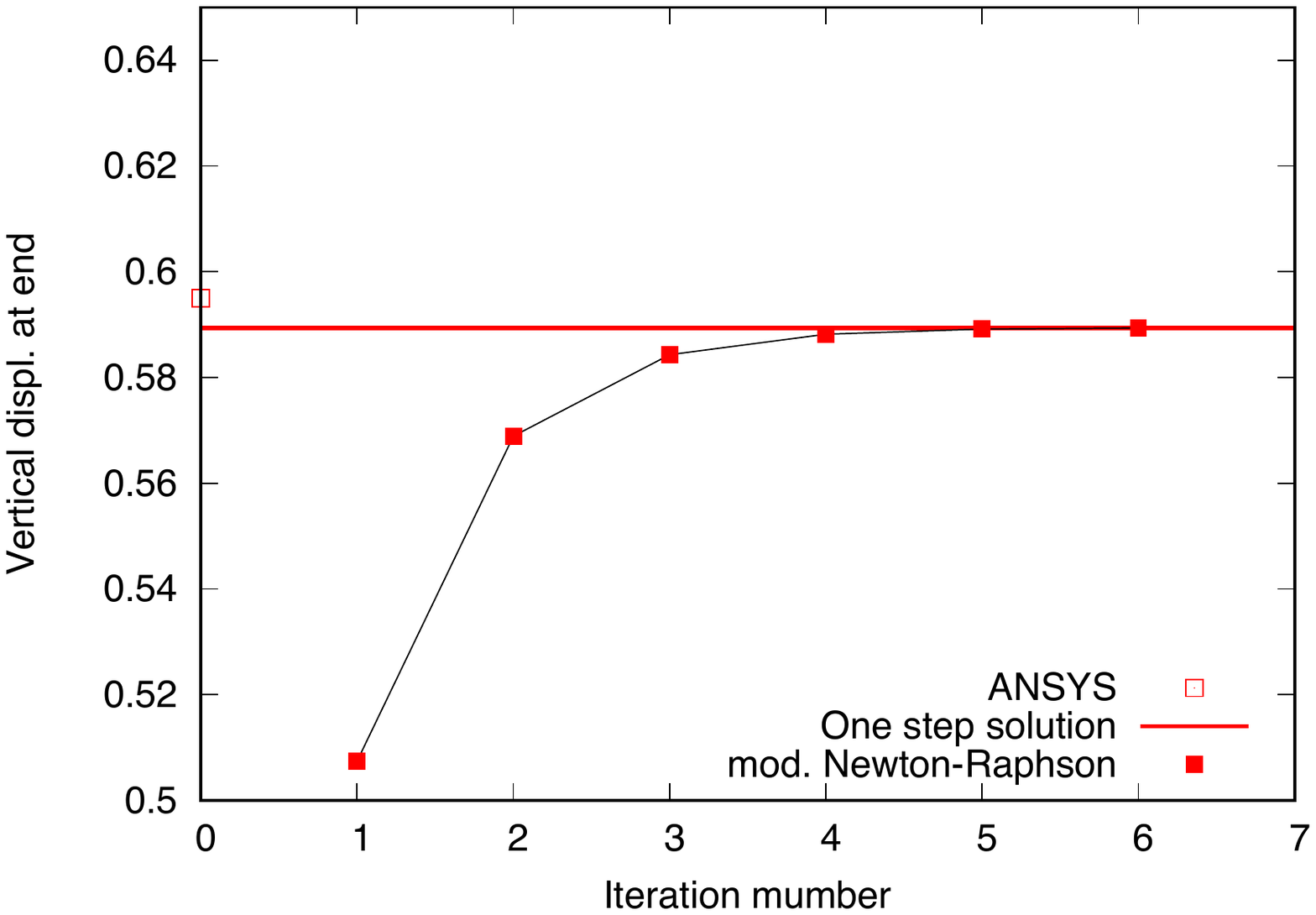}
\caption{Comparison of modified Newton-Raphson solution with one step solution and the result of a FEM analysis}
\label{Ex2conv}
\end{center}
\end{figure}
\subsection{Results}
Only 3 internal points were defined in the horizontal direction for the inclusion. Since the variation of the initial stress is linear in this direction this number was sufficient and for the elastic simulation an increase in the number did not change the result.
The displaced shape of the cantilever is shown in Fig. \ref{Ex2displ}.
In Fig. \ref{Ex2conv} we compare the solution obtained with a modified Newton-Raphson with a one step solution. It can be seen that the iterative solution converges exactly to the one step solution and also agrees well with the result of a FEM analysis with 70323 degrees of freedom using ANSYS.

\section{Summary and Conclusions}
The paper was concerned with the efficient BEM simulation for domains that contain elastic inclusions. Since the BEM relies on fundamental solutions, that exist only for homogeneous and linear domains, this is not a trivial problem. The problem can be solved by considering body forces and this means that volume integrals appear in addition to the surface integrals. The method therefore involves a geometrical discretisation of the volumes or inclusion regions, the computation of strains inside the region and the evaluation of the arising volume integrals.
For the geometrical discretisation of general inclusions we have used already published methods using NURBS bounding surfaces. For the computation of the strains inside the inclusions derived fundamental solutions can be used but because of their high singularity the integration is complicated and not very efficient. This is why we have used a novel approach of computing the strains, by taking the numerical derivatives as is commonly used in the FEM. 

The volume integration can either use body forces, which are derivatives of the initial stresses or the initial stresses directly. To avoid having to take the derivatives the latter strategy has been adopted here. This approach not only makes the evaluation of the volume integral more efficient but also allows a one step solution for the case where the inclusions are elastic and this constitutes a main innovation of this paper.

We propose to use the software in our area of expertise, namely the simulation of underground excavations, where a great number of rock bolts may be used for ground support. It is therefore imperative to implement an efficient integration scheme for bolts. Since numerical integration is only approximate and its' accuracy depends on the number of integration points, we propose to use analytical integration for the bolts, which is fast and accurate. Using reasonable assumptions such as that the cross-sectional area of the bolt is very small compared with the overall dimension of the problem and using local coordinates we can arrive at relatively simple integration results.

Two examples have been presented, one which tests the implementation of the bar inclusion and one that tests the one step solution.
It is shown that the results agree well with comparative solutions obtained by the FEM.

The paper provides a good basis for further work. The next step is to combine the one step solution with non-linear simulations. The left hand side, modified due to the presence of inclusions can now be used, instead of the one for the homogeneous problem, for simulations involving non-linear material behaviour. Indeed, the possibility exists that a true Newton-Raphson method can be applied for non-linear problems, a first for the BEM. 

\bibliographystyle{myplainnat}
\bibliography{bookbib}

\end{document}